\documentclass[11pt,fleqn]{elsarticle}
\usepackage{amssymb}

\usepackage{hyperref}
\usepackage{times}
\usepackage{graphicx}
\usepackage{wrapfig}
\usepackage{fancyhdr}
\usepackage{amsmath}
\usepackage{xcolor}

\newtheorem{theorem}{Theorem}

\newtheorem{lemma}[theorem]{Lemma}

\newenvironment{proof}[1][Proof]{\textbf{#1.} }{\ \rule{0.5em}{0.5em}}
\input{tcilatex}
\begin{document}

\begin{frontmatter}
		
		\title{Recursive Computation of the Multipole Expansions of Layer Potential Integrals over Simplices for Efficient Fast Multipole Accelerated Boundary Elements}
		\author{Nail A. Gumerov\footnote{This author passed away during the submission of this paper.}, Shoken Kaneko and Ramani Duraiswami \\
			Department of Computer Science and  Institute for Advanced Computer Studies\\ University of Maryland, College Park, MD 20742}

		
\begin{abstract}
	In boundary element methods (BEM) in $\mathbb{R}^3$, matrix elements and right hand sides are typically computed via analytical or numerical quadrature of the layer potential multiplied by some function over line, triangle and tetrahedral volume elements. When the problem size gets large, the resulting linear systems are often solved iteratively via Krylov subspace methods, with fast multipole methods (FMM)  used to accelerate the  matrix vector products needed.  When FMM acceleration is used, most  entries of the matrix never need be computed explicitly - {\em they are only  needed in terms of their contribution to the multipole expansion coefficients.}  We propose a new fast method - \emph{Quadrature to Expansion (Q2X)} - for the analytical generation of the multipole expansion coefficients produced  by the integral expressions for single and double layers on surface triangles; charge distributions over line segments and over tetrahedra in the volume; so that the overall method is well integrated into the FMM, with controlled error. The method is based on the $O(1)$ per moment cost recursive computation of the moments. The method is developed for boundary element methods involving the Laplace Green's function in ${\mathbb R}^3$.  The derived recursions are  first compared against classical quadrature algorithms, and then integrated into FMM accelerated boundary element and vortex element methods. Numerical tests are presented and discussed. 

\end{abstract}
	\end{frontmatter}



\section{Introduction}

The fast multipole method (FMM) \cite{Greengard1987:JCP} is often used to accelerate the matrix vector products needed in the Krylov subspace based iterative methods for solution of the large linear systems that result from boundary element methods (BEM) (e.g., for a review see \cite{Nishimura2002:ASME}). The surface is usually discretized in patches, which in the simplest case are plane triangles. For the lowest-order consistent approximation, the functions of interest are approximated via low order polynomials over these triangles, often constant or linear. For constant panel methods, the BEM matrix elements can be represented as integrals of the Green's function and its derivatives over the triangles (``layer potentials''). 
Techniques for evaluating these layer potentials have been extensively developed including analytical methods, e.g., \cite{Newman1986:JEM,Lenoir2012:SISC,Tihon2018:IEEETAP,Wilton1984:IEEETAP,KGD2023RIPE,Gumerov2023:arxiv}. 
In \cite{KGD2023RIPE}, the authors have developed a method for layer potential integral evaluation via dimensionality reduction and recursions for the Laplace and Helmholtz kernel for higher order elements. 

The conventional BEM solution proceeds by evaluating all matrix elements via quadrature over the surface elements and solving the linear system either via matrix decomposition or iteratively. When used with the FMM accelerated iterative methods (FMMBEM), the far field contribution to the matrix vector products are evaluated approximately using expansions that evaluate the product to a specified precision, and the matrix elements are only needed as contributions to the expansion coefficients. While there exist methods to evaluate the element integrals exactly or approximately with controlled accuracy (e.g., \cite{Adelman2016:IEEE}), this  may not be consistent with the accuracy needed in the FMM. We need high precision integrals for some elements (those in the near field) and for the others (in the far field) the integrals are only needed as far as they contribute to the specified accuracy to the consolidated far field multipole expansions. 
For a solver to be efficient, it is important that the accuracy of the computation of all sub-parts of the entire algorithm is kept consistent, so that no redundant or deficient computation occurs. 

Related works on the evaluation of integrals of multipoles include \cite{Nishimura1999:EABE}, which presents a method to analytically evaluate the multipole integrals for the hypersingular Laplace kernel based on dimensionality reduction, which is an approach different from the one based on recursions we will propose in this paper. 
There have also been some alternate proposals for approaches that are numerically more compact. An example of such an approach is in \cite{Adelman2017:IEEE}, where the far field expansions for triangles are simply based on the expansions of equivalent monopoles and dipoles for the evaluation points separated from the source triangle by a sufficient cutoff distance. While this method can be tuned to improve the accuracy of the boundary element method approximation and geometrical errors, it loses the advantage the FMM has of consistent error control in all steps. In this context we also should mention the QBX (quadrature by expansion) method developed in a number of recent studies \cite{Wala_2019-oe, Wala2020-bb, Siegel2018-xj}. Besides establishing the error bounds the developers of the QBX are seeking high-order quadratures, which also can account for higher order shape approximation. This contrasts with the present approach based on low-order shape approximation, analytical quadrature to expansion coefficients, and fast recursive process for generation of multipole expansions.

The focus of this paper is the development of analytical quadrature formulae based on special function recursions that produce the higher order multipole expansion coefficients from the low order ones. These recursions can be used to compute the integrals of multipoles up to the order needed to be consistent with the prescribed FMM error bound, and ensure that the computational and memory needs of the FMMBEM are optimal, and only slightly exceed the cost per iteration of a single evaluation of a FMM matrix vector product for isolated sources located at the panel centers of the BEM mesh. The method is presented for boundary integrals related to the Green's function for the Laplace equation, which  appear in many application areas, including those for linear elasticity \cite{Of2005:CVS}, viscosity \cite{Wang2007:IJNME}, and low frequency electromagnetism \cite{Gumerov2020:IEEE}, and the approach can be used in the FMM for a broad class of problems. 

The BEM in $\mathbb{R}^3$ is used to solve not only the Laplace equation with regular boundaries, but also in domains with filaments and wires (in vortex filament methods and for magnetostatics based on the Biot-Savart law, and for fields around electric wires). Further, when the equations being solved have a volume term, the internal volume needs to be discretized via simplices (tetrahedral elements) and integrals taken over volume elements. We provide formulate for the  multipole coefficients from quadrature for these as well.

\subsection{Summary of results}
We established the following main results for Quadrature to Expansion in this paper:
\begin{itemize}
	\item All integrals for multipole coefficients of layer potential kernels of the Laplace Green's function that arise in boundary element methods based on the collocation method with constant panel approximation are reduced to integrals of harmonic polynomials over unit simplices (line, triangle, tetrahedron).
	\item These integrals can be computed to any order needed by using nested elementary recursions with at most four terms.
	\item Starting exact values for these recursions are provided.
\end{itemize}

Accuracy and performance tests of these recursions were conducted by comparing with numerical and exact quadrature formulae. The recursions were then integrated  into our FMMBEM solver for the Laplace equation and into our vortex element method codes.  Using our approach we obtain a consistently accurate FMMBEM solver whose cost per iteration is very comparable to a simple FMM matrix vector product.

\section{Problem Statement}

\begin{figure}[tbh]
\begin{center}
\includegraphics[width=0.6\textwidth, trim=0.in 0.25in 5.25in
		0.2in]{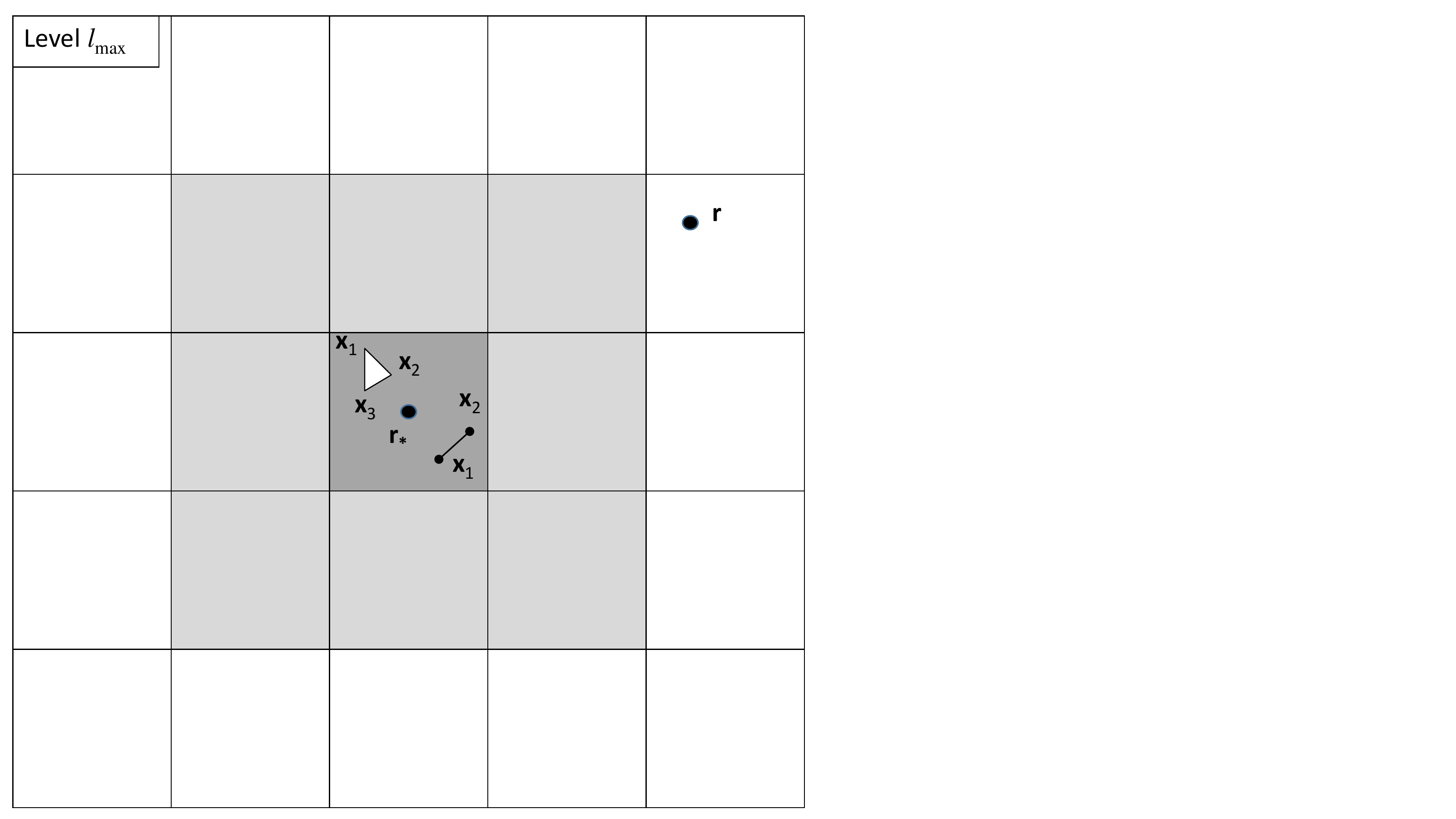}
\end{center}
\caption{An illustrative sketch (in two-dimensions for clarity) of the problem geometry of the use of the FMM for matrix vector products. The finest level of the octree is shown with the source box (the darker shaded box), which contains simplices over which quadrature must be performed (only the line segment and triangle are shown) to achieve the matrix elements. A multipole expansion with origin at the box center needs to be generated for these elements, and should provide accurate representations of the integrals over the simplices to the specified accuracy, for any point $\mathbf{r}$ outside the neighborhood of the source box (the lighter shaded region). }
\label{Fig1}
\end{figure}
The conventional FMM was built for the evaluation of sums of point sources and dipoles distributed in a region. These points are placed into a hierarchical octree data structure (or variants) and which are then used to build valid representations (usually series), with regions of convergence that are determined by the radius of spheres outside or inside of which the truncated expansions converge to a specified precision (e.g., see \cite{Greengard1987:JCP}, \cite{Gumerov2008:JCP}). The minimum box size in these algorithms is chosen to be one that balances the cost of local summation (of points in the finest level box and its neighborhood) with the cost of the expansion and translation costs of the points in the far field. As shown in Fig. \ref{Fig1} we assume that the lines, triangles and tetrahedra over which we will compute integrals in this paper are smaller than the box dimension, and are entirely contained within the box and its neighborhood.

In the BEM, as the surface $S$ is discretized by triangles $S_{j}$, $j=1,...,N$ the single and double layer potentials can be approximated as 
\begin{eqnarray}
L\left( \mathbf{r}\right) &=&\int_{S}\sigma \left( \mathbf{r}^{\prime
}\right) G\left( \mathbf{r},\mathbf{r}^{\prime }\right) dS\left( \mathbf{r}
^{\prime }\right) \approx \sum_{j=1}^{N}\sigma _{j}\int_{S_{j}}G\left( \mathbf{r},\mathbf{r}^{\prime }\right) dS\left( \mathbf{r}^{\prime }\right)=
\sum_{j=1}^{N}\sigma _{j}L_{j}\left( \mathbf{r}\right) ,  \label{st1} \\
M\left( \mathbf{r}\right) &=&\int_{S}\sigma \left( \mathbf{r}^{\prime
}\right) \frac{\partial G\left( \mathbf{r},\mathbf{r}^{\prime }\right) }{\partial n\left( \mathbf{r}^{\prime }\right) }dS\left( \mathbf{r}^{\prime
}\right) \approx \sum_{j=1}^{N}\sigma _{j} \int_{S_{j}}\frac{\partial G\left( \mathbf{r},\mathbf{r}^{\prime }\right) }{\partial n\left( \mathbf{r}^{\prime}\right) } =  \sum_{j=1}^{N}\sigma _{j}  M_{j}\left( \mathbf{r}\right), 
\notag
\end{eqnarray}%
where we used  a constant panel approximation, so $\sigma _{j}$ is the value of the surface density at the $j$th triangle, and $L_{j}\left( \mathbf{r}\right) $ and $M_{j}\left( \mathbf{r}\right) $ are the elementary integrals of the Green's function  $G\left(\mathbf{r},\mathbf{r}^{\prime }\right) $ 
\begin{equation}G\left( \mathbf{r},\mathbf{r}^{\prime }\right)=\frac{1}{4\pi \left| \mathbf{r}-\mathbf{r}^{\prime }\right| }\end{equation}
and its normal derivative over the triangle. For given $\sigma _{j}$ the integrals can be evaluated at any spatial point $\mathbf{r}$ (either in the domain or on the surface). Moreover, the gradients of these potentials may also be needed if we need to compute domain derivatives, or if we are using the Maue identity.

In formulations of magnetostatics and in vortex methods the Biot-Savart integral is needed:
\begin{equation}
\mathbf{H}\left( \mathbf{r}\right) =\frac{1}{4\pi }\int_{C}\frac{\mathbf{I}\left( \mathbf{r}^{\prime }\right) \times \left( \mathbf{r}-\mathbf{r}^{\prime }\right) }{\left| \mathbf{r}-\mathbf{r}^{\prime }\right| ^{3}}dC\left( \mathbf{r}^{\prime }\right) \approx \nabla \times \sum_{j=1}^{N}\mathbf{I}_{j}\int_{C_{j}}G\left( \mathbf{r},\mathbf{r}^{\prime }\right) dl\left( \mathbf{r}^{\prime }\right)  =\nabla \times \sum_{j=1}^{N}\mathbf{I}_{j}K_{j}\left( \mathbf{r}\right) ,  \label{st3} 
\end{equation}
where $C$ is a single contour, or multiple contours, subdivided into $N$ straight line segments (line elements) $C_{j}$, $j=1,...,N$, $\mathbf{I}$ is the current directed along the line element and represented by a constant value $\mathbf{I}_{j}$ for each element. 

Integral equation solution of the Poisson equation with right hand side $f$ requires computation of the integral
\begin{equation}
N\left( \mathbf{r}\right) =\int_{V}f\left( \mathbf{r}^{\prime }\right)
G\left( \mathbf{r},\mathbf{r}^{\prime }\right) dV\left( \mathbf{r}^{\prime
}\right) \approx \sum_{j=1}^{N}f_{j}\int_{V_{j}}G\left(\mathbf{r},\mathbf{r}^{\prime }\right) dV\left( \mathbf{r}^{\prime }\right) = \sum_{j=1}^{N}f_{j} N_{j}\left( \mathbf{r}\right) ,
\label{st3.1} 
\end{equation}%
where we assume that the total volume $V\subset \mathbb{R}^{3}$ is subdivided via tetrahedral elements $V_{j}$ and a piece-wise constant approximation is used for $f\left( \mathbf{r}^{\prime }\right)$.

In comparison to previous FMMBEM implementations  designed for summation of monopoles and dipoles, we  modify the FMMBEM  algorithm to directly generate the multipole expansions that evaluate to these integrals at the location of the box center containing the element (Quadrature to Expansion),  
\begin{equation}
F_{j}\left( \mathbf{r}\right) =\sum_{n=0}^{p-1}\sum_{m=-n}^{n}F_{jn}^{m}\left( \mathbf{r}_{\ast }\right) S_{n}^{m}\left( \mathbf{r}-\mathbf{r}_{\ast }\right) +R_{p}^{(F)},\quad F=K,L,M,N.  \label{st4}
\end{equation}
where $\mathbf{r}_{\ast }$ is the expansion center (the center of the smallest box in the octree containing the center of the $j$th surface, line or volume element, $\mathbf{r}_{j}$, so the expansions are valid at $\left| \mathbf{r}-\mathbf{r}_{\ast }\right| >\left| \mathbf{r}_{j}-\mathbf{r}_{\ast }\right| $), $S_{n}^{m}$ are the spherical basis functions (multipoles) of the Laplace equation, $p$ is the truncation number chosen to provide a given accuracy, $K_{jn}^{m}$, $L_{jn}^{m} $, $M_{jn}^{m}$, and $N_{jn}^{m}$ are the expansion coefficients, and $R_{p}^{(K)}$, $R_{p}^{(L)}$, $R_{p}^{(M)}$, and $R_{p}^{(N)}$ are the residuals due to the truncation of infinite series. We propose a fast recursive process to obtain these coefficients, which results in $O(1)$ complexity per expansion coefficient with a low asymptotic constant. 

Note that this process  can be used also at the evaluation stage for Galerkin method of moments.

\section{Recursions to evaluate Multipole coefficients of the Integrals}

\subsection{Mapping to Unit Simplices}

We parameterize the volume, surface, and line elements considered
using scalar parameters $\left( u,v,w\right) $, $\left( u,v\right) $, and $u$ for the volume, surface, and line elements, respectively. Equations for points in the line, the plane, and  the volume respectively are
\begin{eqnarray}
\mathbf{r} &=&\mathbf{R}\left( u\right) =\mathbf{R}_{u}u+\mathbf{R}%
_{0}\qquad (line),  \label{map1} \\
\quad \mathbf{r} &=&\mathbf{R}\left( u,v\right) =\mathbf{R}_{u}u+\mathbf{R}%
_{v}v+\mathbf{R}_{0}\qquad (plane),  \notag \\
\mathbf{r} &=&\mathbf{R}\left( u,v,w\right) =\mathbf{R}_{u}u+\mathbf{R}_{v}v+%
\mathbf{R}_{w}w+\mathbf{R}_{0}\qquad (space),  \notag
\end{eqnarray}%
where $\mathbf{r}$ is point, while $\mathbf{R}_{u}, \mathbf{R}_{v}, \mathbf{R}_{w}$, and $\mathbf{R}_{0}$ are constant vectors in 
$\mathbb{R}^{3}$. This parametrization also provides a mapping of arbitrary
line segments, triangles, and tetrahedrons on to standard unit simplices: a unit segment with end points at $u=0$ and $u=1$; the  triangle, $T$, with vertices $(0,0)$, $(1,0)$, and $(0,1)$ in the $(u,v)$-plane; and the unit  tetrahedron, $Q$, \ with vertices $(0,0,0)$, $(1,0,0)$, $(0,1,0)$, and $(0,0,1)$ in $\left( u,v,w\right) $-space.

The correspondence of the vertices of the elements provides a unique linear mapping to these unit simplices. These are given by  
\begin{eqnarray}
\mathbf{x}_{1} &=&\mathbf{R}\left( 0\right) =\mathbf{R}_{0},\quad \mathbf{x}%
_{2}=\mathbf{R}\left( 1\right) =\mathbf{R}_{u}+\mathbf{R}_{0}\quad (segment),
\label{map2} \\
\mathbf{x}_{1} &=&\mathbf{R}\left( 0,0\right) =\mathbf{R}_{0},\quad \mathbf{x%
}_{2}=\mathbf{R}\left( 1,0\right) =\mathbf{R}_{u}+\mathbf{R}_{0},\quad 
\mathbf{x}_{3}=\mathbf{R}\left( 0,1\right) =\mathbf{R}_{v}+\mathbf{R}%
_{0}\quad (triangle),  \notag \\
\mathbf{x}_{1} &=&\mathbf{R}\left( 0,0,0\right) =\mathbf{R}_{0},\quad 
\mathbf{x}_{2}=\mathbf{R}\left( 1,0,0\right) =\mathbf{R}_{u}+\mathbf{R}%
_{0}\quad (tetrahedron),  \notag \\
\mathbf{x}_{3} &=&\mathbf{R}\left( 0,1,0\right) =\mathbf{R}_{v}+\mathbf{R}%
_{0},\quad \mathbf{x}_{4}=\mathbf{R}\left( 0,0,1\right) =\mathbf{R}_{w}+%
\mathbf{R}_{0}\text{ \ }(tetrahedron).  \notag
\end{eqnarray}
Given coordinates of the vertices, we can determine%
\begin{equation}
\mathbf{R}_{0}=\mathbf{x}_{1},\quad \mathbf{R}_{u}=\mathbf{x}_{2}-\mathbf{x}%
_{1},\quad \mathbf{R}_{v}=\mathbf{x}_{3}-\mathbf{x}_{1},\quad \mathbf{R}_{w}=%
\mathbf{x}_{4}-\mathbf{x}_{1}  \label{map3}
\end{equation}%
where $\mathbf{R}_{w}$ is defined only for the tetrahedron and $\mathbf{R}%
_{v}$ for the triangle and tetrahedron.
The integrals are then transformed as 
\begin{eqnarray}
L_{j} &=&J^{(j)}\int_{T}G\left( \mathbf{r,R}^{(j)}\left( u,v\right) \right)
dS\left( u,v\right) ,\quad J^{(j)}=\left| \mathbf{R}_{u}^{(j)}\times \mathbf{%
R}_{v}^{(j)}\right| ,  \label{map4} \\
M_{j} &=&J^{(j)}\int_{T}\frac{\partial }{\partial n_{j}}G\left( \mathbf{r},%
\mathbf{R}^{(j)}\left( u,v\right) \right) dS\left( u,v\right) ,\quad  \notag
\\
K_{j} &=&J^{(j)}\int_{0}^{1}G\left( \mathbf{r,R}^{(j)}\left( u\right)
\right) du,\quad J^{(j)}=\left| \mathbf{R}_{u}^{(j)}\right| ,  \notag \\
N_{j} &=&J^{(j)}\int_{Q}G\left( \mathbf{r,R}^{(j)}\left( u,v,w\right)
\right) dV\left( u,v,w\right) ,\quad J^{(j)}=\left| \left( \mathbf{R}%
_{u}^{(j)}\times \mathbf{R}_{v}^{(j)}\right) \cdot \mathbf{R}%
_{w}^{(j)}\right| ,  \notag
\end{eqnarray}%
where $J^{(j)}$ is the Jacobian of transform of the $j$th simplex to the
corresponding unit simplex.

\subsection{Basis functions}
There exist two sorts of spherical basis functions satisfying the Laplace equation, the multipoles - which are singular basis functions at the origin of the reference frame, $S_{n}^{m}\left( \mathbf{r}\right) $, and regular basis functions, or harmonic polynomials, $R_{n}^{m}\left( \mathbf{r}\right) $. These can be expressed via spherical harmonics, and can have different normalization. We use the following definitions:%
\begin{eqnarray}
R_{n}^{m}\left( \mathbf{r}\right) &=&\frac{\left( -1\right) ^{n}i^{\left|
m\right| }}{(n+\left| m\right| )!}r^{n}P_{n}^{\left| m\right| }(\cos \theta
)e^{im\varphi },\quad  \label{bas1} \\
S_{n}^{m}\left( \mathbf{r}\right) &=&i^{-\left| m\right| }(n-\left| m\right|
)!r^{-n-1}P_{n}^{\left| m\right| }(\cos \theta )e^{im\varphi },\quad \quad 
\notag \\
n &=&0,1,2,...,\quad m=-n,...,n,\quad \quad  \notag \\
R_{n}^{m}\left( \mathbf{r}\right) &\equiv &S_{n}^{m}\left( \mathbf{r}\right)
\equiv 0,\quad \left| m\right| >n,  \notag \\
\mathbf{r} &\mathbf{=}&\left( x,y,z\right) =r\left( \sin \theta \cos \varphi
,\sin \theta \sin \varphi ,\cos \theta \right)  \notag
\end{eqnarray}%
where $\left( r,\theta ,\varphi \right) $ are the spherical coordinates of
spatial point $\mathbf{r}$, and $P_{n}^{m}$ are the associated Legendre
functions, defined by the Rodrigues formula (see \cite{Abramowitz1988:Book}%
), 
\begin{equation}
P_{n}^{m}\left( \mu \right) =\frac{\left( -1\right) ^{m}\left( 1-\mu
^{2}\right) ^{m/2}}{2^{n}n!}\frac{d^{m+n}}{d\mu ^{m+n}}\left( \mu
^{2}-1\right) ^{n},\quad n\geqslant 0,\quad m\geqslant 0.  \label{bas2}
\end{equation}%
Note that the basis functions obey the following symmetry,%
\begin{equation}
R_{n}^{-m}\left( \mathbf{r}\right) =\left( -1\right) ^{m}\overline{%
R_{n}^{m}\left( \mathbf{r}\right) },\quad S_{n}^{-m}\left( \mathbf{r}%
\right) =\left( -1\right) ^{m}\overline{S_{n}^{m}\left( \mathbf{r}\right) },
\label{bas2.1}
\end{equation}%
where the bar indicates the complex conjugate. This shows that they need be
computed only for non-negative $m$. 
In these bases the Green's function can be expanded as 
\begin{equation}
G\left( \mathbf{r,r}^{\prime }\right) =\frac{1}{4\pi }\sum_{n=0}^{\infty
}\sum_{m=-n}^{n}(-1)^{n}R_{n}^{-m}\left( \mathbf{r}^{\prime }-\mathbf{r}%
_{\ast }\right) S_{n}^{m}\left( \mathbf{r}-\mathbf{r}_{\ast }\right) ,\quad\left| \mathbf{r}-\mathbf{r}_{\ast }\right| >\left| \mathbf{r}%
^{\prime }-\mathbf{r}_{\ast }\right| .  \label{bas3}
\end{equation}%
Comparing this with equations (\ref{st1})-(\ref{st4}), we can see that%
\begin{eqnarray}
L_{jn}^{m}\left( \mathbf{r}_{\ast }\right) &=&\frac{1}{4\pi }%
(-1)^{n}\int_{S_{j}}R_{n}^{-m}\left( \mathbf{r}^{\prime }-\mathbf{r}_{\ast
}\right) dS\left( \mathbf{r}^{\prime }\right) ,  \label{bas4} \\
M_{jn}^{m}\left( \mathbf{r}_{\ast }\right) &=&\frac{1}{4\pi }(-1)^{n}\mathbf{%
n}_{j}\cdot \int_{S_{j}}\nabla R_{n}^{-m}\left( \mathbf{r}^{\prime }-\mathbf{%
r}_{\ast }\right) dS\left( \mathbf{r}^{\prime }\right) ,  \notag \\
K_{jn}^{m}\left( \mathbf{r}_{\ast }\right) &=&\frac{1}{4\pi }%
(-1)^{n}\int_{C_{j}}R_{n}^{-m}\left( \mathbf{r}^{\prime }-\mathbf{r}_{\ast
}\right) dC\left( \mathbf{r}^{\prime }\right) ,  \notag \\
N_{jn}^{m}\left( \mathbf{r}_{\ast }\right) &=&\frac{1}{4\pi }%
(-1)^{n}\int_{V_{j}}R_{n}^{-m}\left( \mathbf{r}^{\prime }-\mathbf{r}_{\ast
}\right) dV\left( \mathbf{r}^{\prime }\right) ,  \notag
\end{eqnarray}%
where $\mathbf{n}_{j}$ is the normal to the $j$th surface element.

The following derivatives of the harmonic polynomials are also
harmonic polynomials, as the following relations hold (e.g., see \cite%
{Gumerov2006:JCP}), 
\begin{equation}
\frac{\partial }{\partial z}R_{n}^{m}\left( \mathbf{r}\right) =-R_{n-1}^{m}\left( \mathbf{r}\right) , \quad  \frac{\partial }{\partial \eta }R_{n}^{m}\left( \mathbf{r}\right) =iR_{n-1}^{m+1}\left( \mathbf{r}\right) , \quad    \frac{\partial }{\partial \xi }R_{n}^{m}\left( \mathbf{r}\right) =iR_{n-1}^{m-1}\left( \mathbf{r}\right) ,\label{bas5}
\end{equation}%
where%
\begin{equation}
\xi =\frac{x+iy}{2},\quad \eta =\frac{x-iy}{2};\quad x=\xi +\eta ,\quad
y=-i\left( \xi -\eta \right) .  \label{bas6}
\end{equation}%
and%
\begin{equation}
\frac{\partial }{\partial \eta }=\frac{\partial }{\partial x}+i\frac{%
\partial }{\partial y},\quad \frac{\partial }{\partial \xi } = \frac{%
\partial }{\partial x}-i\frac{\partial }{\partial y}.  \label{bas7}
\end{equation}%
Thus, we have%
\begin{equation}
\mathbf{n}\cdot \nabla R_{n}^{m}\left( \mathbf{r}\right) =i \frac{n_{x}}{2}\left[
R_{n-1}^{m+1}\left( \mathbf{r}\right) +R_{n-1}^{m-1}\left( \mathbf{r}\right) %
\right] + \frac{n_{y}}{2} \left[ R_{n-1}^{m+1}\left( \mathbf{r}\right)
-R_{n-1}^{m-1}\left( \mathbf{r}\right) \right] -n_{z}R_{n-1}^{m}\left( 
\mathbf{r}\right) ,  \label{bas8}
\end{equation}%
where $\mathbf{n}=\left( n_{x},n_{y},n_{z}\right) $. This shows that
computation of all the integrals can be reduced to computation of elementary
integrals of the harmonic polynomials
\begin{gather}
a_{jn}^{m}=\int_{Q}c_{jn}^{m}\left( u,v,w\right) dS\left( u,v\right)
=\int_{0}^{1}\int_{0}^{1-u}\int_{0}^{1-u-v}c_{jn}^{m}\left( u,v,w\right)
dwdvdu,\quad   \label{bas9} \\
c_{jn}^{m}\left( u,v,w\right) =R_{n}^{m}\left( \mathbf{R}^{(j)}\left(
u,v,w\right) -\mathbf{r}_{\ast }\right) ,  \notag \\
i_{jn}^{m}=\int_{T}g_{jn}^{m}\left( u,v\right) dS\left( u,v\right)
=\int_{0}^{1}\int_{0}^{1-u}g_{jn}^{m}\left( u,v\right) dvdu,\quad
g_{jn}^{m}\left( u,v\right) =R_{n}^{m}\left( \mathbf{R}^{(j)}\left(
u,v\right) -\mathbf{r}_{\ast }\right) ,  \notag \\
p_{jn}^{m}=\int_{0}^{1}r_{jn}^{m}\left( u\right) du,\quad r_{jn}^{m}\left(
u\right) =R_{n}^{m}\left( \mathbf{R}^{(j)}\left( u\right) -\mathbf{r}_{\ast
}\right) .  \notag
\end{gather}%
The expansion coefficients then can be found as%
\begin{eqnarray}
K_{jn}^{m} &=&\frac{J^{(j)}}{4\pi }(-1)^{n}p_{jn}^{-m},\quad L_{jn}^{m}=%
\frac{J^{(j)}}{4\pi }(-1)^{n}i_{jn}^{-m},\quad   \label{bas10} \\
M_{jn}^{m} &=&\frac{J^{(j)}}{4\pi }(-1)^{n}l_{jn}^{-m},\quad N_{jn}^{m}=%
\frac{J^{(j)}}{4\pi }(-1)^{n}a_{jn}^{-m},  \notag \\
l_{jn}^{-m} &=&i\frac{n_{x}}{2}\left[ i_{j,n-1}^{-m+1}+i_{j,n-1}^{-m-1}\right] + \frac{n_{y}}{2}%
\left[ i_{j,n-1}^{-m+1}-i_{j,n-1}^{-m-1}\right] -n_{z}i_{j,n-1}^{-m}.  \notag
\end{eqnarray}

\subsection{Recursions to evaluate the integrals}
We prove a few lemmas related to the integrals (\ref{bas9}). To simplify notation we drop the subscript/superscript $j$ to denote the  $j$th triangle. Furthermore, we set the origin of the coordinate system so that  $\mathbf{r}_{\ast }=\mathbf{0}.$ This is permissible as, according to Eq. (\ref{map1}), we have%
\begin{eqnarray}
\mathbf{R}\left( u\right) -\mathbf{r}_{\ast } &=&\mathbf{R}_{u}u+\mathbf{R}%
_{0}-\mathbf{r}_{\ast }\qquad (line),  \label{int0} \\
\quad \mathbf{R}\left( u,v\right) -\mathbf{r}_{\ast } &=&\mathbf{R}_{u}u+%
\mathbf{R}_{v}v+\mathbf{R}_{0}-\mathbf{r}_{\ast }\qquad (plane),  \notag \\
\mathbf{R}\left( u,v,w\right) -\mathbf{r}_{\ast } &=&\mathbf{R}_{u}u+\mathbf{%
R}_{v}v+\mathbf{R}_{w}w+\mathbf{R}_{0}-\mathbf{r}_{\ast }\qquad (space), 
\notag
\end{eqnarray}%
This just amounts to shifting the origin $\mathbf{R}_{0}$ in the global coordinate system.

Note that $R_{n}^{m}$ are homogeneous polynomials of degree $n$ of
arguments $\left( x,y,z\right) $ or $\left( \xi ,\eta ,z\right) $. Hence,
according to the Euler's theorem on homogeneous functions we have%
\begin{equation}
nR_{n}^{m}\left( \mathbf{r}\right) =\xi \frac{\partial R_{n}^{m}\left( 
\mathbf{r}\right) }{\partial \xi }+\eta \frac{\partial R_{n}^{m}\left( 
\mathbf{r}\right) }{\partial \eta }+z\frac{\partial R_{n}^{m}\left( \mathbf{r%
}\right) }{\partial z}.  \label{int1}
\end{equation}%
Substituting here differential relations (\ref{bas5}), we obtain%
\begin{equation}
nR_{n}^{m}\left( \mathbf{r}\right) =i\xi R_{n-1}^{m-1}\left( \mathbf{r}%
\right) +i\eta R_{n-1}^{m+1}\left( \mathbf{r}\right) -zR_{n-1}^{m}\left( 
\mathbf{r}\right) .  \label{int2}
\end{equation}

Furthermore, according to Eq. (\ref{map1}) for the points in space, we have%
\begin{eqnarray}
\xi \left( u,v,w\right) &=&\xi _{u}u+\xi _{v}v+\xi _{w}w+\xi _{0},
\label{int3} \\
\eta \left( u,v,w\right) &=&\eta _{u}u+\eta _{v}v+\eta _{w}w+\eta _{0}, 
\notag \\
z\left( u,v,w\right) &=&z_{u}u+z_{v}v+z_{w}w+z_{0},  \notag
\end{eqnarray}%
For points on a surface or a line we have the same relations, where the
dependence on parameter $w$ should be neglected ($\xi _{w}=\eta _{w}=z_{w}=0$%
). Also, for points on the line we have ($\xi _{v}=\eta _{v}=z_{v}=0$).

\begin{lemma}
The following recursion holds for $p_{n}^{m}$ defined by Eq. (\ref{bas9})%
\begin{equation}
p_{n}^{m}=\frac{1}{n+1}\left[ i\xi _{0}p_{n-1}^{m-1}+i\eta
_{0}p_{n-1}^{m+1}-z_{0}p_{n-1}^{m}+q_{n}^{m}\right] ,\quad
q_{n}^{m}=r_{n}^{m}\left( 1\right) .  \label{int5}
\end{equation}
\end{lemma}

\begin{proof}
Due to relations (\ref{bas5}), (\ref{int2}), and (\ref{int3}) at $\xi
_{v}=\eta_{v}=z_{v}= \xi_{w}=\eta_{w}=z_{w}= 0$ we have 
\begin{eqnarray*}
u\frac{\partial r_{n}^{m}}{\partial u} &=&\xi_{u}u\frac{\partial R_{n}^{m}}{\partial \xi } +
\eta_{u}u\frac{\partial R_{n}^{m}}{\partial \eta } + z_{u}u\frac{\partial R_{n}^{m}}{\partial z} \\
&=&\xi \frac{\partial R_{n}^{m}}{\partial \xi }+\eta \frac{\partial R_{n}^{m}%
}{\partial \eta }+z\frac{\partial R_{n}^{m}}{\partial z}-\left( \xi _{0}%
\frac{\partial R_{n}^{m}}{\partial \xi }+\eta _{0}\frac{\partial R_{n}^{m}}{%
\partial \eta }+z_{0}\frac{\partial R_{n}^{m}}{\partial z}\right)  \\
&=&nr_{n}^{m}-\left( i\xi _{0}r_{n-1}^{m-1}+i\eta
_{0}r_{n-1}^{m+1}-z_{0}r_{n-1}^{m}\right) .
\end{eqnarray*}%
This shows that%
\begin{equation*}
\int_{0}^{1}u\frac{\partial r_{n}^{m}}{\partial u}du=np_{n}^{m}-\left( i\xi
_{0}p_{n-1}^{m-1}+i\eta _{0}p_{n-1}^{m+1}-z_{0}p_{n-1}^{m}\right) .
\end{equation*}%
On the other hand, integration by parts results in%
\begin{equation*}
\int_{0}^{1}u\frac{\partial r_{n}^{m}}{\partial u}du=%
\int_{0}^{1}udr_{n}^{m}=r_{n}^{m}\left( 1\right)
-\int_{0}^{1}r_{n}^{m}du=q_{n}^{m}-p_{n}^{m}.
\end{equation*}%
Comparing these expressions we obtain the statement of the lemma.
\end{proof}

\begin{lemma}
The following recursion holds for $i_{n}^{m}$ defined in Eq. (\ref{bas9}),%
\begin{equation}
i_{n}^{m}=\frac{1}{n+2}\left( i\xi _{0}i_{n-1}^{m-1}+i\eta
_{0}i_{n-1}^{m+1}-z_{0}i_{n-1}^{m}+j_{n}^{m}\right) ,  \label{int6}
\end{equation}%
where 
\begin{equation}
j_{n}^{m}=\int_{0}^{1}h_{n}^{m}\left( u\right) du,\quad h_{n}^{m}\left(
u\right) =g_{n}^{m}\left( u,1-u\right) .  \label{int7}
\end{equation}
\end{lemma}

\begin{proof}
Due to relations (\ref{bas5}), (\ref{int2}), and (\ref{int3}) we have 
\begin{eqnarray}
u\frac{\partial g_{n}^{m}}{\partial u}+v\frac{\partial g_{n}^{m}}{\partial v}
&=&\left( \xi _{u}u+\xi _{v}v\right) \frac{\partial R_{n}^{m}}{\partial \xi }%
+\left( \eta _{u}u+\eta _{v}v\right) \frac{\partial R_{n}^{m}}{\partial \eta 
}+\left( z_{u}u+z_{v}v\right) \frac{\partial R_{n}^{m}}{\partial z}
\label{int7.1} \\
&=&\xi \frac{\partial R_{n}^{m}}{\partial \xi }+\eta \frac{\partial R_{n}^{m}%
}{\partial \eta }+z\frac{\partial R_{n}^{m}}{\partial z}-\left( \xi _{0}%
\frac{\partial R_{n}^{m}}{\partial \xi }+\eta _{0}\frac{\partial R_{n}^{m}}{%
\partial \eta }+z_{0}\frac{\partial R_{n}^{m}}{\partial z}\right)   \notag \\
&=&ng_{n}^{m}-\left( i\xi _{0}g_{n-1}^{m-1}+i\eta
_{0}g_{n-1}^{m+1}-z_{0}g_{n-1}^{m}\right) .  \notag
\end{eqnarray}%
This shows that%
\begin{equation}
\int_{T}\left( u\frac{\partial g_{n}^{m}}{\partial u}+v\frac{\partial
g_{n}^{m}}{\partial v}\right) dS=ni_{n}^{m}-\left( i\xi
_{0}i_{n-1}^{m-1}+i\eta _{0}i_{n-1}^{m+1}-z_{0}i_{n-1}^{m}\right) .
\label{int7.2}
\end{equation}%
On the other hand, this integral can be computed as follows%
\begin{eqnarray}
\int_{T}\left( u\frac{\partial g_{n}^{m}}{\partial u}+v\frac{\partial
g_{n}^{m}}{\partial v}\right) dS &=&\int_{0}^{1}\int_{0}^{1-v}u\frac{%
\partial g_{n}^{m}}{\partial u}dudv+\int_{0}^{1}\int_{0}^{1-u}v\frac{%
\partial g_{n}^{m}}{\partial v}dvdu  \label{int7.3} \\
&=&\int_{0}^{1}\int_{0}^{1-v}u\left( d_{u}g_{n}^{m}\right)
dv+\int_{0}^{1}\int_{0}^{1-u}v\left( d_{v}g_{n}^{m}\right) du  \notag \\
&=&\int_{0}^{1}\left[ \left. ug_{n}^{m}\right|
_{0}^{1-v}-\int_{0}^{1-v}g_{n}^{m}du\right] dv+\int_{0}^{1}\left[ \left.
vg_{n}^{m}\right| _{0}^{1-u}-\int_{0}^{1-u}g_{n}^{m}dv\right] du  \notag \\
&=&\int_{0}^{1}\left[ \left( 1-v\right) g_{n}^{m}\left( 1-v,v\right) \right]
dv+\int_{0}^{1}\left[ \left( 1-u\right) g_{n}^{m}\left( u,1-u\right) \right]
du-2\int_{T}g_{n}^{m}dS  \notag \\
&=&\int_{0}^{1}g_{n}^{m}\left( u,1-u\right) du-2i_{n}^{m},  \notag
\end{eqnarray}%
where the last equality can be checked by substitution of variable $v=1-u$
in the integral over $v,$%
\begin{eqnarray}
&&\int_{0}^{1}\left[ \left( 1-v\right) g_{n}^{m}\left( 1-v,v\right) \right]
dv+\int_{0}^{1}\left[ \left( 1-u\right) g_{n}^{m}\left( u,1-u\right) \right]
du  \label{int7.4} \\
&=&\int_{0}^{1}\left[ ug_{n}^{m}\left( u,1-u\right) \right] du+\int_{0}^{1}%
\left[ \left( 1-u\right) g_{n}^{m}\left( u,1-u\right) \right]
du=\int_{0}^{1}g_{n}^{m}\left( u,1-u\right) du.  \notag
\end{eqnarray}%
Equalizing expressions Eq. (\ref{int7.2}) and (\ref{int7.3}) we obtain%
\begin{equation*}
\int_{0}^{1}g_{n}^{m}\left( u,1-u\right)
du-2\int_{T}g_{n}^{m}dS=n\int_{T}g_{n}^{m}dS-\left( i\xi
_{0}\int_{T}g_{n-1}^{m-1}dS+i\eta
_{0}\int_{T}g_{n-1}^{m+1}dS-z_{0}\int_{T}g_{n-1}^{m}dS\right) .
\end{equation*}%
This is equivalent to the statement of the lemma (\ref{int6}).
\end{proof}

\begin{lemma}
The following recursion holds for $j_{n}^{m}$ defined by Eq. (\ref{int7})%
\begin{equation}
j_{n}^{m}=\frac{1}{n+1}\left[ i\left( \xi _{0}+\xi _{v}\right)
j_{n-1}^{m-1}+i\left( \eta _{0}+\eta _{v}\right) j_{n-1}^{m+1}-\left(
z_{0}+z_{v}\right) j_{n-1}^{m}+k_{n}^{m}\right] ,\quad
k_{n}^{m}=g_{n}^{m}\left( 1,0\right) .  \label{int8}
\end{equation}
\end{lemma}

\begin{proof}
Note that at $\mathbf{r}_{\ast }=\mathbf{0}$ the integrand in Eq. (\ref{int7}%
) is 
\begin{eqnarray}
h_{n}^{m}\left( u\right) &=&g_{n}^{m}\left( u,1-u\right) =R_{n}^{m}\left( 
\mathbf{R}\left( u,1-u\right) \right)  \label{int8.1} \\
&=&R_{n}^{m}\left( \mathbf{R}_{u}u+\mathbf{R}_{v}\left( 1-u\right) +\mathbf{R%
}_{0}\right) =R_{n}^{m}\left( \left( \mathbf{R}_{u}-\mathbf{R}_{v}\right) u+%
\mathbf{R}_{v}+\mathbf{R}_{0}\right) .  \notag
\end{eqnarray}%
This is the same as $r_{n}^{m}\left( u\right) $ for the line integral, but
with $\mathbf{R}_{u}-\mathbf{R}_{v}$ instead of $\mathbf{R}_u$ and $\mathbf{R}_{v}+%
\mathbf{R}_{0}$ instead of $\mathbf{R}_{0}$. From Lemma 1 (Eq. \ref{int5})
we obtain the statement of the lemma by replacing $\xi _{0},\eta _{0}$, and $%
z_{0}$ with $\xi _{0}+\xi _{v},\eta _{0}+\eta _{v},$ and $z_{0}+z_{v}$. We
also note that $q_{n}^{m}$ there corresponds to $k_{n}^{m}$ since at $u=1$
we have $g_{n}^{m}\left( 1,0\right) =r_{n}^{m}\left( 1\right)
=R_{n}^{m}\left( \mathbf{R}_{u}+\mathbf{R}_{0}\right) .$
\end{proof}

\begin{lemma}
The following recursion holds for $a_{n}^{m}$ defined in Eq. (\ref{bas9}),%
\begin{equation}
a_{n}^{m}=\frac{1}{n+3}\left( i\xi _{0}a_{n-1}^{m-1}+i\eta
_{0}a_{n-1}^{m+1}-z_{0}a_{n-1}^{m}+b_{n}^{m}\right) ,  \label{int9}
\end{equation}%
where 
\begin{equation}
b_{n}^{m}=\int_{0}^{1}\int_{0}^{1-v}d_{n}^{m}\left( u,v\right) dudv,\quad
d_{n}^{m}\left( u,v\right) =c_{n}^{m}\left( u,v,1-u-v\right) .  \label{int10}
\end{equation}
\end{lemma}

\begin{proof}
Due to relations (\ref{bas5}), (\ref{int2}), and (\ref{int3}) we have 
\begin{eqnarray}
&&u\frac{\partial c_{n}^{m}}{\partial u}+v\frac{\partial c_{n}^{m}}{\partial
v}+w\frac{\partial c_{n}^{m}}{\partial w}  \label{int10.1} \\
&=&\left( \xi _{u}u+\xi _{v}v+\xi _{w}w\right) \frac{\partial R_{n}^{m}}{%
\partial \xi }+\left( \eta _{u}u+\eta _{v}v+\eta _{w}w\right) \frac{\partial
R_{n}^{m}}{\partial \eta }+\left( z_{u}u+z_{v}v+z_{w}w\right) \frac{\partial
R_{n}^{m}}{\partial z}  \notag \\
&=&\xi \frac{\partial R_{n}^{m}}{\partial \xi }+\eta \frac{\partial R_{n}^{m}%
}{\partial \eta }+z\frac{\partial R_{n}^{m}}{\partial z}-\left( \xi _{0}%
\frac{\partial R_{n}^{m}}{\partial \xi }+\eta _{0}\frac{\partial R_{n}^{m}}{%
\partial \eta }+z_{0}\frac{\partial R_{n}^{m}}{\partial z}\right)   \notag \\
&=&nc_{n}^{m}-\left( i\xi _{0}c_{n-1}^{m-1}+i\eta
_{0}c_{n-1}^{m+1}-z_{0}c_{n-1}^{m}\right) .  \notag
\end{eqnarray}%
So, 
\begin{equation}
\int_{Q}\left( u\frac{\partial c_{n}^{m}}{\partial u}+v\frac{\partial
c_{n}^{m}}{\partial v}+w\frac{\partial c_{n}^{m}}{\partial w}\right)
dV=na_{n}^{m}-\left( i\xi _{0}a_{n-1}^{m-1}+i\eta
_{0}a_{n-1}^{m+1}-z_{0}a_{n-1}^{m}\right) .  \label{int10.2}
\end{equation}%
Integral of the left hand side of (\ref{int10.2}) over the standard
tetrahedron can be represented as 
\begin{eqnarray*}
&&\int_{Q}\left( u\frac{\partial c_{n}^{m}}{\partial u}+v\frac{\partial
c_{n}^{m}}{\partial v}+w\frac{\partial c_{n}^{m}}{\partial w}\right) dV \\
&=&\int_{0}^{1}\int_{0}^{1-v}\int_{0}^{1-v-w}u\frac{\partial c_{n}^{m}}{%
\partial u}dudwdv+\int_{0}^{1}\int_{0}^{1-w}\int_{0}^{1-u-w}v\frac{\partial
c_{n}^{m}}{\partial v}dvdudw+\int_{0}^{1}\int_{0}^{1-u}\int_{0}^{1-u-v}w%
\frac{\partial c_{n}^{m}}{\partial w}dwdvdu \\
&=&\int_{0}^{1}\int_{0}^{1-v}\left[ \left. uc_{n}^{m}\right|
_{0}^{1-v-w}-\int_{0}^{1-v-w}c_{n}^{m}du\right] dwdv+\int_{0}^{1}%
\int_{0}^{1-w}\left[ \left. vc_{n}^{m}\right|
_{0}^{1-u-w}-\int_{0}^{1-u-w}c_{n}^{m}dv\right] dudw \\
&&+\int_{0}^{1}\int_{0}^{1-v}\left[ \left. wc_{n}^{m}\right|
_{0}^{1-u-v}-\int_{0}^{1-u-v}c_{n}^{m}dw\right] dudv \\
&=&\int_{0}^{1}\int_{0}^{1-v}\left[ \left( 1-v-w\right) c_{n}^{m}\left(
1-v-w,v,w\right) \right] dwdv+\int_{0}^{1}\int_{0}^{1-w}\left[ \left(
1-u-w\right) c_{n}^{m}\left( u,1-u-w,w\right) \right] dudw \\
&&+\int_{0}^{1}\int_{0}^{1-v}\left[ \left( 1-u-v\right) c_{n}^{m}\left(
u,v,1-u-v\right) \right] dvdu-3\int_{Q}c_{n}^{m}dV \\
&=&\int_{0}^{1}\int_{0}^{1-v}c_{n}^{m}\left( u,v,1-u-v\right)
dudv-3\int_{Q}c_{n}^{m}dV.
\end{eqnarray*}%
where the last identity can be checked by substitution $w=1-u-v$ for
integrals over $w$ (at fixed $v$ for the first and at fixed $u$ for the
second integral). Equalizing the last expression with Eq. (\ref{int10.2}) we
obtain the statement of the lemma.
\end{proof}

\begin{lemma}
The following recursion holds for $b_{n}^{m}$ defined by Eq. (\ref{int10})%
\begin{eqnarray}
b_{n}^{m} &=&\frac{1}{n+2}\left( i\left( \xi _{0}+\xi _{w}\right)
b_{n-1}^{m-1}+i\left( \eta _{0}+\eta _{w}\right) b_{n-1}^{m+1}-\left(
z_{0}+z_{w}\right) b_{n-1}^{m}+e_{n}^{m}\right) ,  \label{int11} \\
e_{n}^{m} &=&\frac{1}{n+1}\left[ i\left( \xi _{0}+\xi _{v}\right)
e_{n-1}^{m-1}+i\left( \eta _{0}+\eta _{v}\right) e_{n-1}^{m+1}-\left(
z_{0}+z_{v}\right) e_{n-1}^{m}+f_{n}^{m}\right] ,\quad
f_{n}^{m}=c_{n}^{m}\left( 1,0,0\right) .  \notag
\end{eqnarray}
\end{lemma}

\begin{proof}
Note that at $\mathbf{r}_{\ast }=\mathbf{0}$ the integrand in Eq. (\ref%
{int10}) is 
\begin{eqnarray}
d_{n}^{m}\left( u,v\right)  &=&c_{n}^{m}\left( u,v,1-u-v\right)
=R_{n}^{m}\left( \mathbf{R}\left( u,v,1-u-v\right) \right)   \label{int12} \\
&=&R_{n}^{m}\left( \mathbf{R}_{u}u+\mathbf{R}_{v}v+\mathbf{R}_{w}\left(
1-u-v\right) +\mathbf{R}_{0}\right)   \notag \\
&=&R_{n}^{m}\left( \left( \mathbf{R}_{u}-\mathbf{R}_{w}\right) u+\left( 
\mathbf{R}_{v}-\mathbf{R}_{w}\right) v+\mathbf{R}_{w}+\mathbf{R}_{0}\right) .
\notag
\end{eqnarray}%
This is the same as $g_{n}^{m}\left( u,v\right) $ for the surface integral,
but with $\mathbf{R}_{u}-\mathbf{R}_{w}$, $\mathbf{R}_{v}-\mathbf{R}_{w}$,
and $\mathbf{R}_{0}+\mathbf{R}_{w}$ instead of $\mathbf{R}_{u}$, $\mathbf{R}%
_{v}$, and $\mathbf{R}_{0}$. Hence, the statement of this lemma is a
corollary of Lemma 2 and Lemma 3 since $f_{n}^{m}=c_{n}^{m}\left(
1,0,0\right) =g_{n}^{m}\left( 1,0\right) =r_{n}^{m}\left( 1\right)
=R_{n}^{m}\left( \mathbf{R}_{u}+\mathbf{R}_{0}\right) $.
\end{proof}

\subsection{Recursions}

Lemmas 1-5 along with relation (\ref{int2}) deliver a recursive method to
compute all integrals $p_{n}^{m}$ , $i_{n}^{m}$, and $a_{n}^{m}$. 
Auxiliary integrals $j_{n}^{m},$ $e_{n}^{m},$ and $b_{n}^{m}$ should also be
computed recursively. Separate recursions are needed for $q_{n}^{m}$, $%
k_{n}^{m}$, and $f_{n}^{m}$. Since 
\begin{eqnarray}
q_{n}^{m} &=&r_{n}^{m}\left( 1\right) =g_{n}^{m}\left( 1,0\right)
=c_{n}^{m}\left( 1,0,0\right) =k_{n}^{m}=f_{n}^{m}  \label{rec1} \\
&=&R_{n}^{m}\left( \xi _{u}+\xi _{0},\eta _{u}+\eta _{0},z_{u}+z_{0}\right) 
\notag
\end{eqnarray}%
we obtain from Eq. (\ref{int2})
\begin{equation}
q_{n}^{m}=\frac{1}{n}\left[ i\left( \xi _{0}+\xi _{u}\right)
q_{n-1}^{m-1}+i\left( \eta _{0}+\eta _{u}\right) q_{n-1}^{m+1}-\left(
z_{0}+z_{u}\right) q_{n-1}^{m}\right] .  \label{rec2}
\end{equation}%
We note also that $j_{n}^{m}$ and $e_{n}^{m}$ satisfy the same recursions
with the same initial values. So, these quantities are simply the same.

Note that in all recursions $a_{n}^{m},b_{n}^{m}$,$i_{n}^{m}$, $j_{n}^{m}$,
and $q_{n}^{m}$ should be set to zero for $\left| m\right| >n$. This simply
follows from definition of these quantities and the fact that $%
R_{n}^{m}\left( \mathbf{r}\right) \equiv 0$ for $\left| m\right| >n$. This
means that all functions of degree $n$ can be obtained from functions of
degree $n-1$. The starting values for the recursions can be computed
directly, since $R_{0}^{0}\equiv 1$, and so%
\begin{eqnarray}
q_{0}^{0} &=&R_{0}^{0}=1,  \label{rec3} \\
p_{0}^{0} &=&j_{0}^{0}=\int_{0}^{1}R_{0}^{0}du=1,  \notag \\
i_{0}^{0} &=&b_{0}^{0}=\int_{0}^{1}\int_{0}^{1-u}R_{0}^{0}dvdu=\frac{1}{2}, 
\notag \\
a_{0}^{0} &=&\int_{0}^{1}\int_{0}^{1-u}\int_{0}^{1-u-v}R_{0}^{0}dvdudw=\frac{%
1}{6}.  \notag
\end{eqnarray}%
Summarizing, we obtained the following recursions for the line integrals,%
\begin{eqnarray}
q_{n}^{m} &=&\frac{1}{n}\left[ i\left( \xi _{0}+\xi _{u}\right)
q_{n-1}^{m-1}+i\left( \eta _{0}+\eta _{u}\right) q_{n-1}^{m+1}-\left(
z_{0}+z_{u}\right) q_{n-1}^{m}\right] ,  \label{rec4} \\
p_{n}^{m} &=&\frac{1}{n+1}\left[ i\xi _{0}p_{n-1}^{m-1}+i\eta
_{0}p_{n-1}^{m+1}-z_{0}p_{n-1}^{m}+q_{n}^{m}\right] ,\quad n=1,2,...  \notag
\end{eqnarray}%
The surface integrals can be computed using the recursions%
\begin{eqnarray}
q_{n}^{m} &=&\frac{1}{n}\left[ i\left( \xi _{0}+\xi _{u}\right)
q_{n-1}^{m-1}+i\left( \eta _{0}+\eta _{u}\right) q_{n-1}^{m+1}-\left(
z_{0}+z_{u}\right) q_{n-1}^{m}\right] ,  \label{rec5} \\
j_{n}^{m} &=&\frac{1}{n+1}\left[ i\left( \xi _{0}+\xi _{v}\right)
j_{n-1}^{m-1}+i\left( \eta _{0}+\eta _{v}\right) j_{n-1}^{m+1}-\left(
z_{0}+z_{v}\right) j_{n-1}^{m}+q_{n}^{m}\right] ,  \notag \\
i_{n}^{m} &=&\frac{1}{n+2}\left( i\xi _{0}i_{n-1}^{m-1}+i\eta
_{0}i_{n-1}^{m+1}-z_{0}i_{n-1}^{m}+j_{n}^{m}\right) ,\quad n=1,2,...  \notag
\end{eqnarray}%
For the volume integrals we have%
\begin{eqnarray}
q_{n}^{m} &=&\frac{1}{n}\left[ i\left( \xi _{0}+\xi _{u}\right)
q_{n-1}^{m-1}+i\left( \eta _{0}+\eta _{u}\right) q_{n-1}^{m+1}-\left(
z_{0}+z_{u}\right) q_{n-1}^{m}\right] ,  \label{rec5.1} \\
j_{n}^{m} &=&\frac{1}{n+1}\left[ i\left( \xi _{0}+\xi _{v}\right)
j_{n-1}^{m-1}+i\left( \eta _{0}+\eta _{v}\right) j_{n-1}^{m+1}-\left(
z_{0}+z_{v}\right) j_{n-1}^{m}+q_{n}^{m}\right] ,  \notag \\
b_{n}^{m} &=&\frac{1}{n+2}\left[ i\left( \xi _{0}+\xi _{w}\right)
b_{n-1}^{m-1}+i\left( \eta _{0}+\eta _{w}\right) b_{n-1}^{m+1}-\left(
z_{0}+z_{w}\right) b_{n-1}^{m}+j_{n}^{m}\right] ,  \notag \\
\quad a_{n}^{m} &=&\frac{1}{n+3}\left( i\xi _{0}a_{n-1}^{m-1}+i\eta
_{0}a_{n-1}^{m+1}-z_{0}a_{n-1}^{m}+b_{n}^{m}\right) ,\quad n=1,2,...  \notag
\end{eqnarray}

Note then, that some other bases and symmetries can be exploited to make
computations more efficient. For example, in \cite{Gumerov2008:JCP} the
entire FMM was implemented on GPUs that lacked complex number implementation using real basis functions defined as%
\begin{equation}
\underline{R_{n}^{m}}=\left\{ 
\begin{array}{c}
\func{Re}\left\{ \widetilde{R}_{n}^{m}\right\} ,\quad m\geqslant 0, \\ 
\func{Im}\left\{ \widetilde{R}_{n}^{m}\right\} ,\quad m<0%
\end{array}%
\right. ,\quad \underline{S_{n}^{m}}=\left\{ 
\begin{array}{c}
\func{Re}\left\{ \widetilde{S}_{n}^{m}\right\} ,\quad m\geqslant 0, \\ 
\func{Im}\left\{ \widetilde{S}_{n}^{m}\right\} ,\quad m<0%
\end{array}%
\right. ,  \label{rec6}
\end{equation}%
where the complex valued basis functions with tildes are related to
functions (\ref{bas1}) via%
\begin{equation}
\widetilde{R}_{n}^{m}\left( \mathbf{r}\right) =i^{\left| m\right|
}R_{n}^{m}\left( \mathbf{r}\right) ,\quad \widetilde{S}_{n}^{m}\left( 
\mathbf{r}\right) =i^{-\left| m\right| }S_{n}^{m}\left( \mathbf{r}\right) .
\label{rec7}
\end{equation}%
The functions obey the following symmetry%
\begin{equation}
\widetilde{R}_{n}^{-m}\left( \mathbf{r}\right) =\overline{\widetilde{R}%
_{n}^{m}\left( \mathbf{r}\right) },\quad \widetilde{S}_{n}^{-m}\left( 
\mathbf{r}\right) =\overline{\widetilde{S}_{n}^{m}\left( \mathbf{r}\right) }.
\label{rec8}
\end{equation}%
Certainly, integration of functions preserves the symmetry. Also integrals
decorated similarly to the integrand satisfy relations (\ref{rec6}) and (\ref%
{rec7}), so for the real basis we have%
\begin{equation}
\underline{F_{n}^{m}}=\left\{ 
\begin{array}{c}
\func{Re}\left\{ \widetilde{F}_{n}^{m}\right\} ,\quad m\geqslant 0, \\ 
\func{Im}\left\{ \widetilde{F}_{n}^{m}\right\} ,\quad m<0%
\end{array}%
\right. =\left\{ 
\begin{array}{c}
\func{Re}\left\{ \widetilde{F}_{n}^{m}\right\} ,\quad m\geqslant 0, \\ 
-\func{Im}\left\{ \widetilde{F}_{n}^{-m}\right\} ,\quad m<0%
\end{array}%
\right. ,\quad F=p,q,i,j,a,b.  \label{rec9}
\end{equation}%
This shows that the complex valued functions need be computed only for
non-negative $m$.

We have then for the line integrals,%
\begin{eqnarray}
\widetilde{q}_{n}^{m} &=&\frac{1}{n}\left[ -\left( \xi _{0}+\xi _{u}\right) 
\widetilde{q}_{n-1}^{m-1}+\left( \eta _{0}+\eta _{u}\right) \widetilde{q}%
_{n-1}^{m+1}-\left( z_{0}+z_{u}\right) \widetilde{q}_{n-1}^{m}\right] ,
\label{rec10} \\
\widetilde{p}_{n}^{m} &=&\frac{1}{n+1}\left[ -\xi _{0}\widetilde{p}%
_{n-1}^{m-1}+\eta _{0}\widetilde{p}_{n-1}^{m+1}-z_{0}\widetilde{p}_{n-1}^{m}+%
\widetilde{q}_{n}^{m}\right] ,\quad m>0,\quad n=1,2,...  \notag
\end{eqnarray}%
\begin{eqnarray}
\widetilde{q}_{n}^{0} &=&\frac{1}{n}\left[ 2\func{Re}\left\{ \left( \eta
_{0}+\eta _{u}\right) \widetilde{q}_{n-1}^{1}\right\} -\left(
z_{0}+z_{u}\right) \widetilde{q}_{n-1}^{0}\right] ,  \label{rec11} \\
\widetilde{p}_{n}^{0} &=&\frac{1}{n+1}\left[ 2\func{Re}\left\{ \eta _{0}%
\widetilde{p}_{n-1}^{1}\right\} -z_{0}\widetilde{p}_{n-1}^{0}+\widetilde{q}%
_{n}^{0}\right] ,\quad n=1,2,..,.  \notag
\end{eqnarray}%
where for relations at $m=0$ we used symmetry (\ref{rec8}) and the fact that 
$\overline{\xi }=\eta $.

The surface integrals can be computed using the same recursions for $%
\widetilde{q}_{n}^{m}$ and the following relations for the other quantities%
\begin{eqnarray}
\widetilde{j}_{n}^{m} &=&\frac{1}{n+1}\left[ -\left( \xi _{0}+\xi
_{v}\right) \widetilde{j}_{n-1}^{m-1}+\left( \eta _{0}+\eta _{v}\right) 
\widetilde{j}_{n-1}^{m+1}-\left( z_{0}+z_{v}\right) \widetilde{j}_{n-1}^{m}+%
\widetilde{q}_{n}^{m}\right] ,  \label{rec12} \\
\widetilde{i}_{n}^{m} &=&\frac{1}{n+2}\left( -\xi _{0}\widetilde{i}%
_{n-1}^{m-1}+\eta _{0}\widetilde{i}_{n-1}^{m+1}-z_{0}\widetilde{i}_{n-1}^{m}+%
\widetilde{j}_{n}^{m}\right) ,\quad m>0,\quad n=1,2,...  \notag
\end{eqnarray}%
\begin{eqnarray}
\widetilde{j}_{n}^{0} &=&\frac{1}{n+1}\left[ 2\func{Re}\left\{ \left( \eta
_{0}+\eta _{v}\right) \widetilde{j}_{n-1}^{1}\right\} -\left(
z_{0}+z_{v}\right) \widetilde{j}_{n-1}^{0}+\widetilde{q}_{n}^{0}\right] ,
\label{rec14} \\
\widetilde{i}_{n}^{0} &=&\frac{1}{n+2}\left( 2\func{Re}\left\{ \eta _{0}%
\widetilde{i}_{n-1}^{1}\right\} -z_{0}\widetilde{i}_{n-1}^{0}+\widetilde{j}%
_{n}^{0}\right) ,\quad n=1,2,...  \notag
\end{eqnarray}

The volume integrals can be computed using the same relations for $%
\widetilde{q}_{n}^{m}$ and $\widetilde{j}_{n}^{m}$ and the following
relations for the other quantities%
\begin{eqnarray}
\widetilde{b}_{n}^{m} &=&\frac{1}{n+2}\left[ -\left( \xi _{0}+\xi
_{w}\right) \widetilde{b}_{n-1}^{m-1}+\left( \eta _{0}+\eta _{w}\right) 
\widetilde{b}_{n-1}^{m+1}-\left( z_{0}+z_{w}\right) \widetilde{b}_{n-1}^{m}+%
\widetilde{j}_{n}^{m}\right] ,  \label{rec14.1} \\
\widetilde{a}_{n}^{m} &=&\frac{1}{n+3}\left( -\xi _{0}\widetilde{a}%
_{n-1}^{m-1}+\eta _{0}\widetilde{a}_{n-1}^{m+1}-z_{0}\widetilde{a}_{n-1}^{m}+%
\widetilde{b}_{n}^{m}\right) ,\quad m>0,\quad n=1,2,...  \notag
\end{eqnarray}%
\begin{eqnarray}
\widetilde{b}_{n}^{0} &=&\frac{1}{n+2}\left[ 2\func{Re}\left\{ \left( \eta
_{0}+\eta _{w}\right) \widetilde{b}_{n-1}^{1}\right\} -\left(
z_{0}+z_{w}\right) \widetilde{b}_{n-1}^{0}+\widetilde{j}_{n}^{0}\right] ,
\label{rec14.2} \\
\widetilde{a}_{n}^{0} &=&\frac{1}{n+3}\left( 2\func{Re}\left\{ \eta _{0}%
\widetilde{a}_{n-1}^{1}\right\} -z_{0}\widetilde{a}_{n-1}^{0}+\widetilde{b}%
_{n}^{0}\right) ,\quad n=1,2,...  \notag
\end{eqnarray}

Since $\widetilde{F}_{0}^{0}=F_{0}^{0},$ $F=p,q,i,j,a,b$, the initial values
for the recursions are provided by Eq. (\ref{rec3}).

Finally, we note that integrals for the double layer potential (\ref{bas10})
can be computed with the recursive process for the single layer potential
using relations%
\begin{eqnarray}
\widetilde{l}_{n}^{m} &=&\frac{1}{2}\widetilde{i}_{n-1}^{m+1}(n_{x}-in_{y})-%
\frac{1}{2}\widetilde{i}_{n-1}^{m-1}(n_{x}+in_{y})-\widetilde{i}%
_{n-1}^{m}n_{z},\quad m>0,\quad n=1,2,...  \label{rec15} \\
\widetilde{l}_{n}^{0} &=&\func{Re}\left[ \widetilde{i}%
_{n-1}^{1}(n_{x}-in_{y})\right] -\widetilde{i}_{n-1}^{0}n_{z},\quad n=1,2,...
\notag
\end{eqnarray}
\section{Numerical tests}

We conducted several numerical tests to confirm the validity of the obtained
formulae and measure the performance.

\subsection{Accuracy tests }

Integrals $L,M,$ and $K$ in equations (\ref{st1}) and (\ref{st3}) can be
computed analytically. In the numerical tests we compared the results
obtained using analytical expressions and the expansions truncated at
different truncation numbers. In all tests the expansion center is fixed at
the origin of the reference frame, while the evaluation point is set at the
point with coordinates%
\begin{equation}
\mathbf{r}=d\left( \frac{\sqrt{3}}{2}, 0, \frac{1}{2}\right) ,  \label{num0}
\end{equation}%
where the distance between the expansion center and the evaluation point, $d$, is varying. The center of the simplex is placed at $\mathbf{r}_{c}=\left( 
\sqrt{3}/2,0,0\right) $, while the simplex vertices are located
on a sphere of radius $r_{t}:$ 
\begin{eqnarray}
\mathbf{x}_{1} &=&\mathbf{r}_{c}+r_{t}\left( -1,0,0\right) ,\quad \mathbf{x}%
_{2}=\mathbf{r}_{c}+r_{t}\left( 1,0,0\right) \quad (segment),  \label{num0.1}
\\
\mathbf{x}_{1} &=&\mathbf{r}_{c}+r_{t}\left( 1,0,0\right) ,\quad \mathbf{x}%
_{2,3}=\mathbf{r}_{c}+r_{t}\left( -\frac{1}{2},\pm \frac{\sqrt{3}}{2}%
,0\right) \quad (triangle),  \notag \\
\mathbf{x}_{1} &=&\mathbf{r}_{c}+r_{t}\left( 1,0,0\right) ,\quad \mathbf{x}%
_{2,3}=\mathbf{r}_{c}+r_{t}\left( -\frac{1}{3},-\frac{\sqrt{2}}{3},\pm \sqrt{%
\frac{2}{3}}\right) ,\quad (tetrahedron),  \notag \\
\mathbf{x}_{4} &=&\mathbf{r}_{c}+r_{t}\left( -\frac{1}{3},\frac{2\sqrt{2}}{3}%
,0\right) \text{ }(tetrahedron),  \notag
\end{eqnarray}%
These locations correspond to the place where the FMM multipole expansions have maximum asymptotic error, using the octree data structure, in which case one also should use $d=1.5$. We did some accuracy tests by varying $r_{t}$ in range $r_{t}\in $[0.01, 0.2], but found that the effect of varying of the truncation number $p$ and $d$ is much stronger. Below
we presented results for $r_{t}=0.1$, which is typical for FMMBEM simulations.

\subsubsection{Segment}
For the line integral we have
\begin{eqnarray}
K &=&\frac{1}{4\pi }\int_{\mathbf{x}_{1}}^{\mathbf{x}_{2}}\frac{dC\left( 
\mathbf{r}^{\prime }\right) }{\left| \mathbf{r}-\mathbf{r}^{\prime }\right| }
=\frac{1}{4\pi }\int_{-1}^{1}\frac{d\xi }{\sqrt{\xi ^{2}-4\zeta \xi \cos
\alpha +4\zeta ^{2}}}  \label{num1}  \\
&=&\frac{1}{4\pi }\ln \frac{\sqrt{4\zeta ^{2}-4\zeta \cos \alpha +1}%
+1-2\zeta \cos \alpha }{\sqrt{4\zeta ^{2}+4\zeta \cos \alpha +1}-1-2\zeta
\cos \alpha }.  \notag
\end{eqnarray}%
Here we performed the following transformation of the original integral%
\begin{eqnarray}
\mathbf{r}^{\prime } &=&\mathbf{r}_{0}+\frac{1}{2}\mathbf{p}\xi ,\quad 
\mathbf{r}_{0}=\frac{1}{2}\left( \mathbf{x}_{1}+\mathbf{x}_{2}\right) ,\quad 
\mathbf{p}=\mathbf{x}_{2}-\mathbf{x}_{1},  \label{num2} \\
\quad \xi &\in &\left[ -1,1\right] ,\quad dl=\frac{1}{2}J\xi ,\quad J=\left| 
\mathbf{p}\right| ,  \notag
\end{eqnarray}%
and denoted%
\begin{equation}
\zeta =\frac{\left| \mathbf{r}-\mathbf{r}_{0}\right| }{J},\quad \cos \alpha =%
\frac{\left( \mathbf{r}-\mathbf{r}_{0}\right) \cdot \mathbf{p}}{\left| 
\mathbf{r}-\mathbf{r}_{0}\right| J}.  \label{num3}
\end{equation}

\begin{figure}[tbh]
\begin{center}
\includegraphics[width=0.96\textwidth, trim=1in 1in 0.25in
		0.24in]{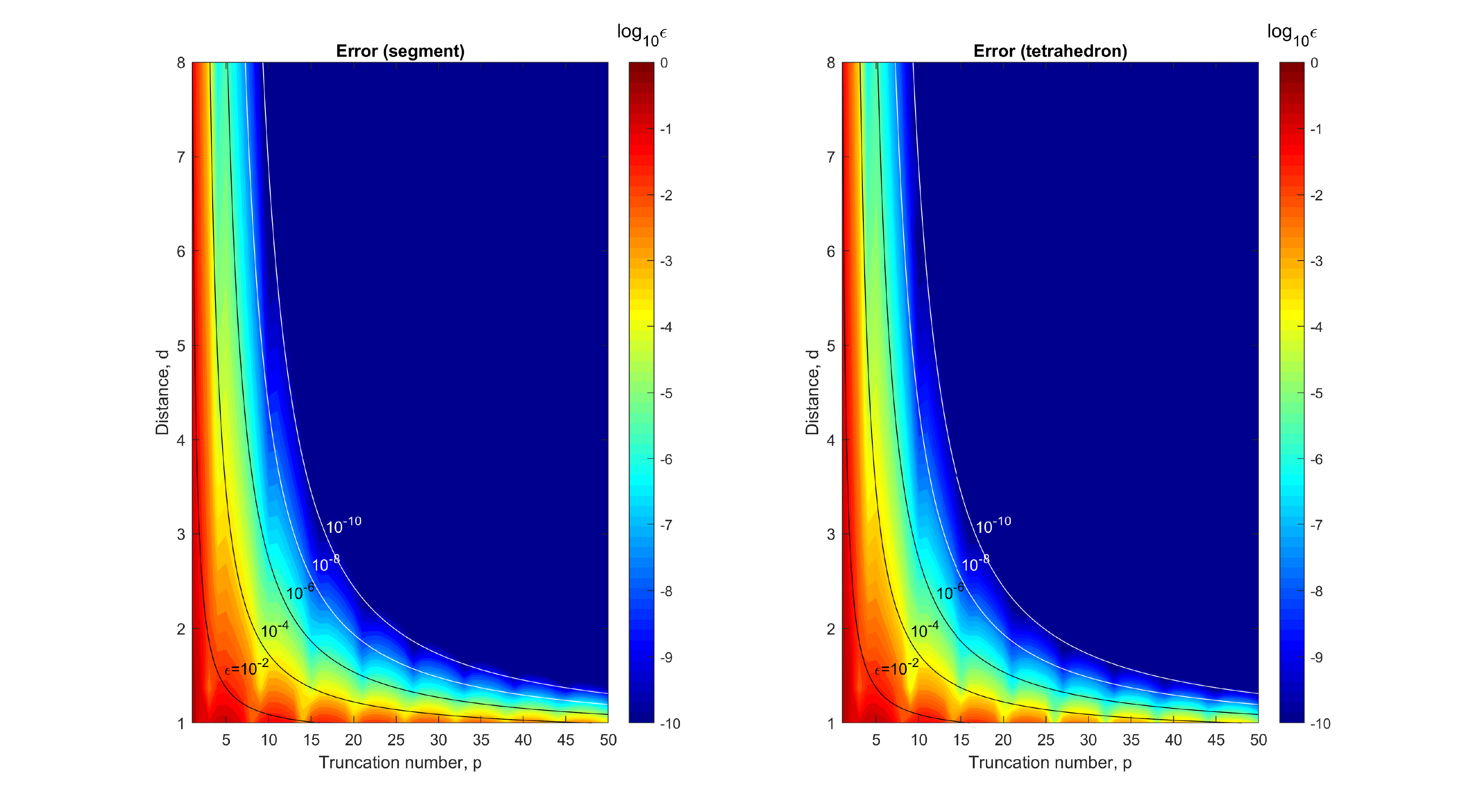}
\end{center}
\caption{Computed errors in the evaluation of the line and tetrahedron
integrals at distance $d$ through the multipole expansions derived for
various values of the truncation number $p$. The lines show the theoretical
error levels according to Eq. (\ref{num4}). The oscillatory pattern is a property of the multipole expansion and is not unique to the proposed method. }
\label{Fig2}
\end{figure}

Figure \ref{Fig2} shows the relative error, $\epsilon =\left|
K-K_{an}\right| /K_{an}$, where $K_{an}$ is provided by equations (\ref{num1}%
)-(\ref{num3}) and $K$ is computed using the series expansion and the
recursive procedure (\ref{rec10})-(\ref{rec11}) for the expansion
coefficients. It is seen that the error can be reduced to a desired level by
a proper selection of parameter $p$. The distance $d$ plays also an
important role. Theoretically, the error depends on $p$ and $d$ as%
\begin{equation}
\epsilon =C\left( \frac{\left| \mathbf{r}_{c}-\mathbf{r}_{\ast }\right| }{d}%
\right) ^{p},  \label{num4}
\end{equation}%
where $C$ is some constant, since the series can be majorated by a geometric
progression$.$ We plotted these curves on Fig. \ref{Fig2} at $C=0.1$ and
different fixed values of $\epsilon $. It is seen that the computed error is
consistent with the theoretical error bound.

\subsubsection{Triangle}

While analytical expressions for integrals for flat quadrilaterals can be
found in \cite{Hess1967:PAS} and can be modified for triangles, we can
obtain analytical expressions by reduction of the surface integrals to the
contour integrals using the Gauss divergence theorem and then by evaluating
the integrals over the segments using the primitives. Below we display only
the final expressions \cite{Adelman2016-ik} needed for validation of the
obtained recursions. 
\begin{eqnarray}
L\left( \mathbf{r}\right) &=&\int_{S}G\left( \mathbf{r},\mathbf{r}^{\prime
}\right) dS\left( \mathbf{r}^{\prime }\right) =\frac{1}{4\pi }\sum_{q=1}^{3}%
\left[ L_{P}\left( l_{q}-x_{q},h,z_{q}\right) -L_{P}\left(
-x_{q},h,z_{q}\right) \right] ,  \label{num5} \\
M\left( \mathbf{r}\right) &=&\int_{S}\frac{\partial G\left( \mathbf{r},%
\mathbf{r}^{\prime }\right) }{\partial n\left( \mathbf{r}^{\prime }\right) }%
dS\left( \mathbf{r}^{\prime }\right) =\frac{1}{4\pi }\sum_{q=1}^{3}\left[
M_{P}\left( l_{q}-x_{q},h,z_{q}\right) -M_{P}\left( -x_{q},h,z_{q}\right) %
\right] ,  \notag \\
x_{q} &=&\left( \mathbf{r-x}_{q}\right) \cdot \mathbf{i}_{q},\quad h=\left|
\left( \mathbf{r-x}_{1}\right) \cdot \mathbf{n}\right| ,\quad z_{q}=\left( 
\mathbf{r-x}_{q}\right) \cdot \mathbf{n}_{q},  \notag \\
\mathbf{i}_{q} &=&\frac{1}{l_{q}}\left( \mathbf{r}_{q(\func{mod}3)+1}-%
\mathbf{x}_{q}\right) ,\quad l_{q}=\left| \mathbf{x}_{q(\func{mod}N_{e})+1}-%
\mathbf{x}_{q}\right| ,\quad \mathbf{n}_{q}=\mathbf{i}_{q}\times \mathbf{n},
\notag \\
\quad L_{P}\left( x,y,z\right) &=&-yM_{P}\left( x,y,z\right) -z\ln \left|
r+x\right| ,\quad M_{P}\left( x,y,z\right) =-\arctan \frac{xz(r-y)}{%
rz^{2}+yx^{2}},  \notag \\
r &=&\sqrt{x^{2}+y^{2}+z^{2}}.  \notag
\end{eqnarray}
\begin{figure}[tbh]
\begin{center}
\includegraphics[width=0.96\textwidth, trim=1in .5in 1in
		0.24in]{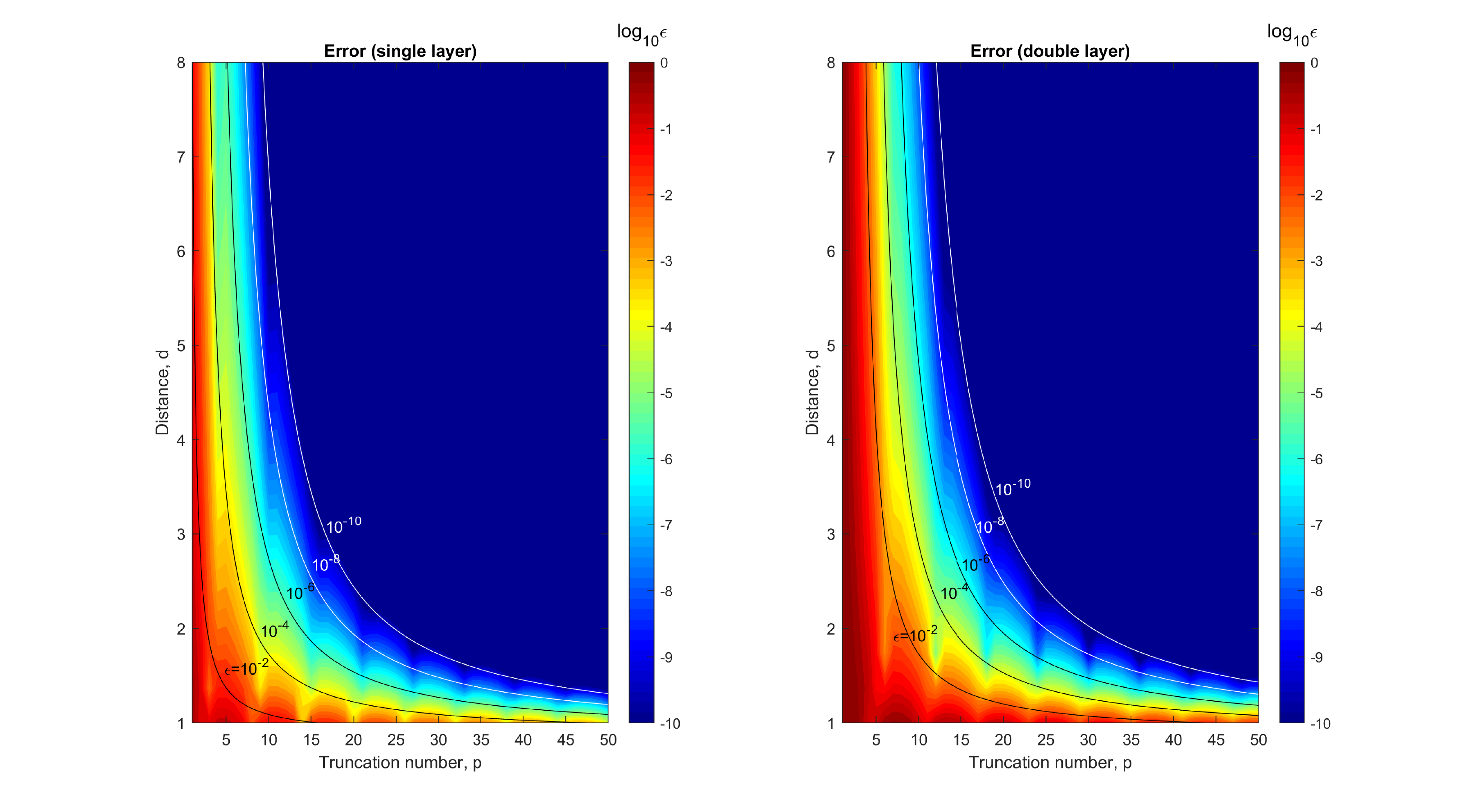}
\end{center}
\caption{Error in the evaluation of the single layer and double layer
integrals over a boundary triangle at distance $d$ through the multipole
expansions derived for various values of the truncation number $p$. The
lines show the theoretical error levels according to Eqs (\ref{num4}) and (%
\ref{num6}).}
\label{Fig3}
\end{figure}

Figure \ref{Fig3} shows the relative errors between the analytical solutions
and $L$ and $M$ computed using the series expansion and the recursive
procedure (\ref{rec12})-(\ref{rec15}) for the expansion coefficients. Note
that while the error for the single layer potential can be described by Eq. (%
\ref{num4}), the error of the double layer potential can be estimated as%
\begin{equation}
\epsilon =Cp\left( \frac{\left| \mathbf{r}_{c}-\mathbf{r}_{\ast }\right| }{d}%
\right) ^{p-1},  \label{num6}
\end{equation}%
with the same constant $C$ as for the single layer. This relation is
obtained by estimation of the residual of the derivative of the geometric
progression majorating the single layer potential$.$ We plotted curves
described by Eqs (\ref{num4}) and (\ref{num6}) on Fig. \ref{Fig3} at $C=0.1$
and different fixed values of $\epsilon $. It is seen that the computed
error is consistent with the theoretical error bound.

\subsubsection{Tetrahedron}

Analytical expression for tetrahedrons, can be obtained from those for flat
triangles. Indeed, we have according to the Gauss divergence theorem

\begin{eqnarray}
N\left( \mathbf{r}\right) &=&\int_{V}G\left( \mathbf{r},\mathbf{r}^{\prime
}\right) dV\left( \mathbf{r}^{\prime }\right) =-\frac{1}{2}\int_{V}\nabla _{%
\mathbf{r}^{\prime }}\cdot \left[ \left( \mathbf{r}-\mathbf{r}^{\prime
}\right) G\left( \mathbf{r},\mathbf{r}^{\prime }\right) \right] dV\left( 
\mathbf{r}^{\prime }\right)  \label{num7} \\
&=&-\frac{1}{2}\int_{S}\mathbf{n}\left( \mathbf{r}^{\prime }\right) \cdot %
\left[ \left( \mathbf{r}-\mathbf{r}^{\prime }\right) G\left( \mathbf{r},%
\mathbf{r}^{\prime }\right) \right] dS\left( \mathbf{r}^{\prime }\right) =-%
\frac{1}{2}\sum_{j=1}^{4}\int_{S_{j}}\mathbf{n}_{j}\cdot \left( \mathbf{r}-%
\mathbf{r}^{\prime }\right) G\left( \mathbf{r},\mathbf{r}^{\prime }\right)
dS\left( \mathbf{r}^{\prime }\right) ,  \notag
\end{eqnarray}%
where $S_{j}$ are the faces of tetrahedron $V$ and $\mathbf{n}_{j}$ are the
normals to the faces directed outside the tetrahedron. Denoting $\mathbf{r}%
_{cj}$ the centers of the faces, and noticing that $\mathbf{n}_{j}\cdot
\left( \mathbf{r}^{\prime }-\mathbf{r}_{cj}\right) =0$, we obtain%
\begin{eqnarray}
N\left( \mathbf{r}\right) &=&-\frac{1}{2}\sum_{j=1}^{4}\int_{S_{j}}\mathbf{n}%
_{j}\cdot \left( \mathbf{r-r}_{cj}-\left( \mathbf{r}^{\prime }-\mathbf{r}%
_{cj}\right) \right) G\left( \mathbf{r},\mathbf{r}^{\prime }\right) dS\left( 
\mathbf{r}^{\prime }\right)  \label{num8} \\
&=&-\frac{1}{2}\sum_{j=1}^{4}\mathbf{n}_{j}\cdot \left( \mathbf{r-r}%
_{cj}\right) \int_{S_{j}}G\left( \mathbf{r},\mathbf{r}^{\prime }\right)
dS\left( \mathbf{r}^{\prime }\right) =-\frac{1}{2}\sum_{j=1}^{4}\mathbf{n}%
_{j}\cdot \left( \mathbf{r-r}_{cj}\right) L_{j}\left( \mathbf{r}\right) . 
\notag
\end{eqnarray}%
Here the single layer potential $L_{j}\left( \mathbf{r}\right) $ can be
computed using Eq. (\ref{num5}).

Figure \ref{Fig2} shows the errors for the integral over the tetrahedron.
The theoretical error bound here is provided by Eq. (\ref{num4}) and plotted
in the figure for $C=0.1$ and different fixed values of $\epsilon $. It is
seen that the computed error is consistent with the theoretical error bound.
Also, they agree well with the errors for the line segment, which has the
same bound.

\subsection{Performance}

\subsubsection{Gauss quadrature}

Expansion coefficients $K_{n}^{m},L_{n}^{m},$ and $M_{n}^{m}$ can be
computed exactly using quadrature formulae, such as the Gauss-Legendre
quadrature. Indeed the basis functions $R_{n}^{m}$ are polynomials of degree 
$n$, so, it is sufficient to use the Gauss quadrature of the order of $%
(n+1)/2$ (take the ceiling if this number is not an integer). For triangles
special quadrature formulae of this type can be derived. However, as soon as
the integral can be transformed to an integral over a standard triangles,
which can be represented by a double integral, the Gauss quadrature can be
applied to each integral.

\begin{figure}[tbh]
\begin{center}
\includegraphics[width=0.96\textwidth, trim=0.95in .5in 0.95in
		0.24in]{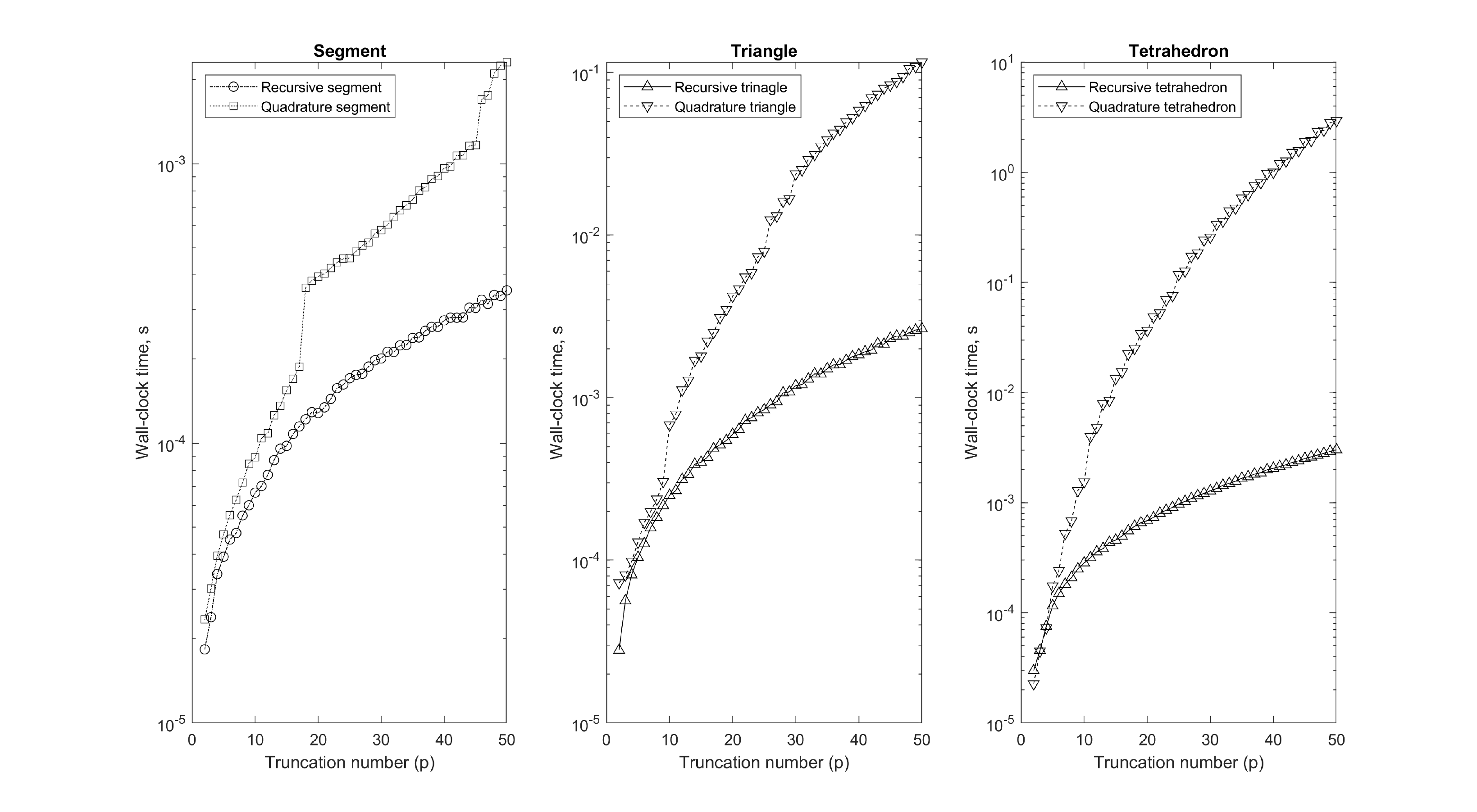}
\end{center}
\caption{The wall-clock times (for Matlab function prototypes) for
computation of the multipole expansions for uniform charge distributions
over a line segment, triangle, and tetrahedron. Two methods are compared for
each case: the present method (recursive) and the method based on the Gauss
quadrature with the number of nodes providing exact integration of the
polynomials of degree $p$. The difference between the results of the two
methods are within the round-off errors for double precision arithmetic. A
significant speedup is seen.}
\label{Fig4}
\end{figure}

First, we compared the results obtained using the present recursive processes and those obtained using the numerical quadrature and confirmed that the difference is of the order of the round-off errors for double precision computing. Second, we measured the time to generate multipole expansion using the quadrature formulae and the recursions. These tests were performed using Matlab, where the quadrature results were computed using vectorization, while the recursions were obtained in a serial manner (loop)
for  increasing $n$. For quadrature, all $R_{n}^{m}$ at the quadrature
nodes were computed recursively from $n=0$ to a given $p-1$ using recursions
for $k_{n}^{m}$ (see Eqs (\ref{rec1}) and (\ref{rec2})) and then summation
was performed using precomputed quadrature weights. Even though the quadrature formulae are computed using vectorization, Fig. \ref{Fig4} shows that the present recursive process outperforms the method based on quadrature. This is clearly seen for the segment, where the wall clock time can be several time smaller for the recursive method. The effect is much stronger for the triangle, where the present process can be 1-2 orders of magnitude faster than the quadrature computations. It is even stronger for the tetrahedron when the recursive computations can be 3 orders of magnitude faster compared to the quadrature method. Of course, the acceleration of computations depends on the truncation number, and for small truncation numbers the benefit is smaller. However, in all cases the  recursive method is observed to be faster.

\subsubsection{Fast multipole method}
While tests in Matlab are fine for demonstrating accuracy, the real test would be in a production code, that we use in our research in fluid flow and computational electromagnetics \cite{Gumerov2008:JCP}. We considered two example computations using the FMM. The times were obtained on a workstation equipped with Intel Xeon CPU processors at 2.1 GHz with total 32 physical cores and 128 GB RAM using an Open MP parallelized code without GPU as described in \cite{Gumerov2008:JCP}%
. 
\begin{figure}[tbh]
\begin{center}
\includegraphics[width=6in]{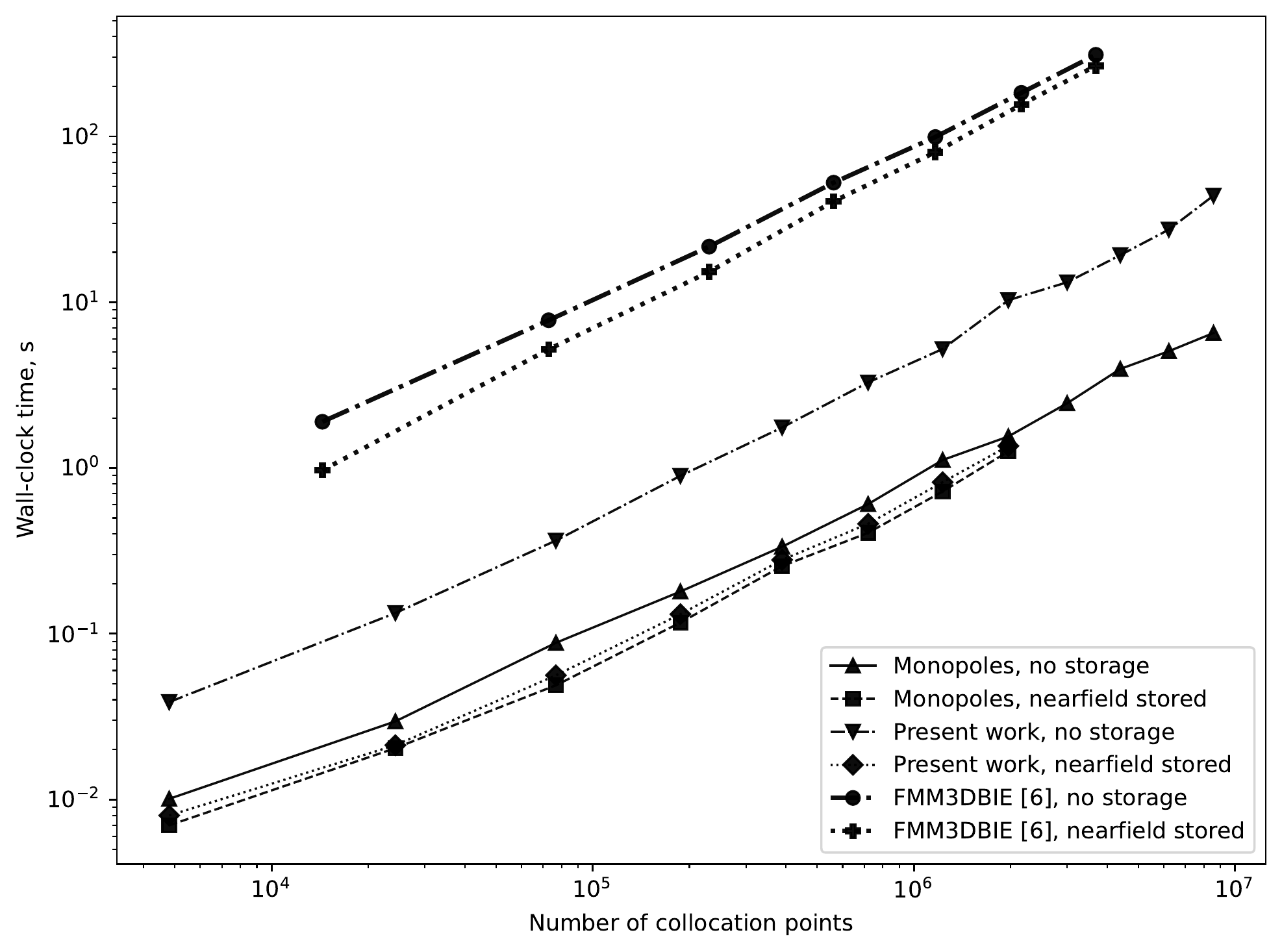}
\end{center}
\caption{Comparisons of the wall-clock time performance of the FMM used for summation of monopoles distributed at the flat triangle centers over the surface of some object and summation of single layer integrals for the same surface covered by triangles. Performance for two options for the present code, and for FMM3DBIE~\cite{Greengard2021:JCP} are shown. The two options are with and without storage of the nearfield interaction sums. The difference in computation times for the option with storage relates only to the speed of generation of multipole expansions. Times are measured for a 32 core Intel Xeon CPU with a clock speed of 2.1 GHz and 128 GB RAM. The penalty for accurate consistent quadrature in comparison with a pure monopole FMM is seen to be minimal for the present work. 
}
\label{Fig5}
\end{figure}

In the first example, we computed the single layer potential $L$ for a unit
sphere with a random distribution of the charge density. The recursive
computations of the multipole expansions were implemented and the
performance was compared with summation of monopoles for the same source and
receiver locations (at the triangle centers). The reason for such test is to
see how much the FMM slows down, due to the computation of integrals. Note
here, that the FMM can be optimized in terms of performance by appropriate
selection of the clustering parameter$,s$, (which is the maximum number of
sources in a box at the maximum level of the octree; this parameters
controls the depth of the octree) (e.g., see \cite{Gumerov2008:JCP}) and we
report only data for the optimized cases. We also note that the matrix for
nearfield interactions (a sparse matrix) in the FMM can be computed ``on the
fly'' every time when the matrix-vector product is needed or it can be
computed and stored. The latter method is preferable for iterative solutions
of large systems, as only the input vector is changing in such a process,
while the matrix does not change. The former method is preferable when only
one or very few evaluations are needed (e.g., for evaluation of the
integrals with already known charge density at some field points). Sometimes
one is forced to use this method even for iterative computations if the 
available RAM is insufficient to store the nearfield interaction matrix. 
In the solver used in this experiment, also The nearfield integrals are evaluated analytically using a dimensionality reduction approach~\cite{Adelman2016-ik}. 

Accordingly, we considered two cases for comparisons shown in Fig. \ref%
{Fig5}. In all cases the truncation number was $p=10$, which provides the
overall relative $L_{2}$-error of the FMM of the order of $\epsilon _{2}\sim
10^{-6}$ (see also \cite{Gumerov2008:JCP}). In all cases the error was the
same for the monopole summation and for summation of the integrals over
triangles. ``On the fly'' computing of the nearfield interaction shows that
the runtime for the monopole summation ($s=300$) is about 3-6 times smaller
than for summation of integrals ($s=30$). However, for the case when the
nearfield matrices are precomputed and stored the monopole summation ($s=400$%
) is only 3-17\% faster than that for the integrals ($s=400$). Note that in
this case the entire difference in performance of the code for monopoles and
triangles is due to the multipole expansions, as all other steps of the
algorithm are the same. A big difference in the performance when the
nearfield interactions are computed directly is only due to direct
computation of integrals (\ref{num5}), which is much more arithmetically
intensive than computation of the Green function.

To compare against a state of the art BEM solver, we computed the same single layer potential using 
FMM3DBIE~\cite{Greengard2021:JCP}. The wall-clock time of the execution was measured for 
the subroutines ``lpcomp\_lap\_comb\_dir\_addsub'' and ``getnearquad\_lap\_comb\_dir'', 
which correspond to the farfield and nearfield evaluations, respectively.
In Fig. \ref{Fig5}, ``no storage" represents the sum of the execution time of a single call to these two functions and ``nearfield stored" represents the execution time of a single call to the farfield evaluation function.
First order elements with three Vioreanu-Rokhlin nodes~\cite{Vioreanu2014} per element, 
which are not shared with neighboring elements, were used. This results in three times more collocation nodes in total, 
compared to the center panel collocation used in the previous experiment. 
The error control parameter $\epsilon$ was set to $10^{-6}$. 
The relative $L_{2}$-error of the computed single layer potential evaluated at the collocation points with respect to the analytical result obtained by the expression (\ref{num5}) was $2.1\times10^{-6}$ and $9.0\times10^{-6}$ for the proposed recursive method and the FMM3DBIE code, respectively. 
As the FMM-BEM is a highly complicated algorithm consisting of many building blocks, the performance depends on various design choices and implementation details of the subroutines constituting the entire algorithm. Also, different solvers provide different functionalities. For example, FMM3DBIE supports curved elements while our solver does not. The performance difference across the different solvers compared here should therefore be considered only as a reference. 

\begin{figure}[tbh]
\begin{center}
\includegraphics[width=0.96\textwidth, trim=2in .5in 2in 0.24in]{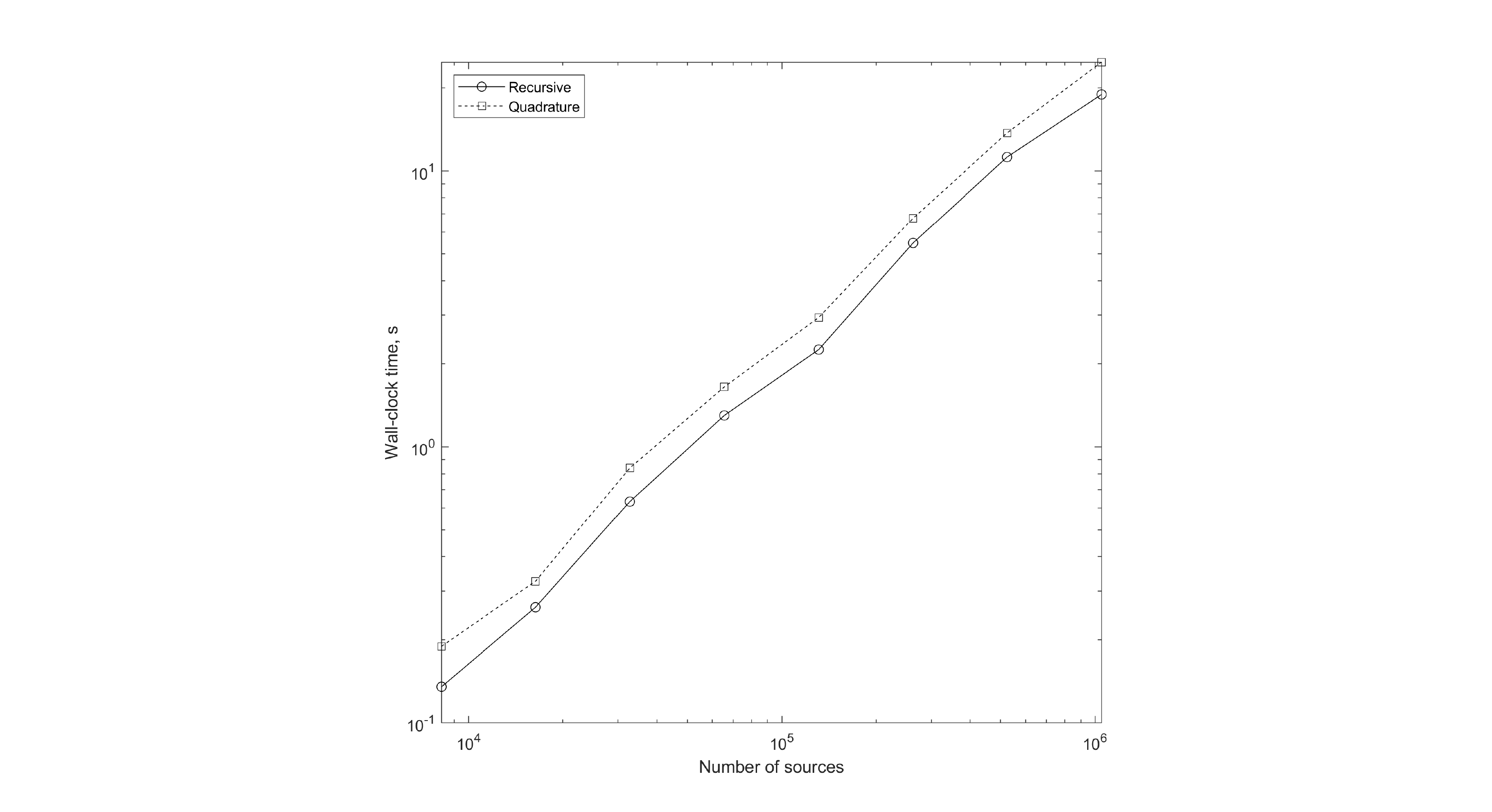}
\end{center}
\caption{The wall-clock times for the FMM accelerated vortex filament method
(truncation number $p=34$, clustering parameter, $s=600$, just summation)
with recursive computation of multipole expansions and computations of the
same expansion using the Gauss quadrature. Times are measured for the same
workstation as in Fig. \ref{Fig5}. A consistent significant speedup is
observed.}
\label{Fig6}
\end{figure}
In the second example, we modified the FMM-accelerated vortex filament
method (see \cite{Gumerov2013:JCP}) with a recursive process described in
the present paper. The wall clock time was measured for different number of
elementary vortex filaments randomly distributed in space and compared with
the times obtained when the multipole expansions of line integrals are
computed using the Gauss quadrature. Figure \ref{Fig6} shows the results for
truncation number $p=34$ and $s=600$. In all cases the use of the recursive
integral computations showed acceleration, which for the results plotted
varied in the range 20-40\%.

\section{Conclusions}
A method for recursive computation of the multipole expansions of the integrals of the single and double potentials over the triangles, of  potentials over line segments, and tetrahedrons is developed and tested. We call this method Quadrature to Expansion (Q2X). The tests show consistency of the results with analytical expressions and with numerical integration using quadrature formulae of appropriate order. In all cases considered the developed process shows significantly better performance in terms of speed and memory than other methods.

According to our tests using an FMM accelerated matrix vector product, the performance of one iteration of the FMMBEM method (summation of integrals) based on the present method of generation of multipole expansions is only 3-17\% slower than summation of multipoles in the case when the nearfield interaction matrix is precomputed and stored. This is a remarkably small overhead. The use of the recursive process may increase the performance of the vortex filament codes up to 40\%.

Note that the present method can be also used in computations of the Galerkin integrals (e.g., see \cite{Adelman2016:IEEE}). Indeed, the no-touch cases can be computed using the multipole expansions. Furthermore, if the Galerkin method is used with the FMM, the far field interactions are locally represented by series over functions $R_{n}^{m}\left( \mathbf{r}\right) $, which should be integrated over the receiving triangles. In this case we have exactly the case considered in the present paper and the fast recursive evaluation of the integrals can be applied here as well.

\section{Acknowledgments}

This work is supported by Cooperative Research Agreements W911NF1420118 and  W911NF2020213 between the University of Maryland and the Army Research Laboratory, with David Hull, Ross Adelman and Steven Vinci as Technical monitors. 
Shoken Kaneko thanks Japan Student Services Organization and Watanabe Foundation for scholarships.

\begin{thebibliography}{19}
	
	\bibitem[\protect\citename{Abramowitz {\em et~al.}, }1988]{Abramowitz1988:Book}
	{\sc Abramowitz, Milton, Stegun, Irene~A, \& Romer, Robert~H}. 1988.
	\newblock Handbook of Mathematical Functions with Formulas, Graphs, and
	Mathematical Tables.
	\newblock {\em Am. J. Phys.}, {\bf 56}(10), 958--958.
	
	\bibitem[\protect\citename{Adelman {\em et~al.}, }2016]{Adelman2016:IEEE}
	{\sc Adelman, R, Gumerov, N~A, \& Duraiswami, R}. 2016.
	\newblock Computation of Galerkin double surface integrals in the {3-D}
	boundary element method.
	\newblock {\em IEEE Transactions on Antennas and Propagation},  {\bf 64}(6), 2389--2400.
	
	\bibitem[\protect\citename{Adelman {\em et~al.}, }2017]{Adelman2017:IEEE}
	{\sc Adelman, R, Gumerov, N~A, \& Duraiswami, R}. 2017.
	\newblock {FMM/GPU-Accelerated} Boundary Element Method for Computational
	Magnetics and Electrostatics.
	\newblock {\em IEEE Trans. Magn.}, {\bf 53}(12), 1--11.
	
	\bibitem[\protect\citename{Adelman, }2016]{Adelman2016-ik}
	{\sc Adelman, Ross}. 2016.
	\newblock {\em Fast and accurate boundary element methods in three dimensions}.
	\newblock Ph.D. thesis, University of Maryland, College Park.
	
	\bibitem[\protect\citename{Greengard \& Rokhlin, }1987]{Greengard1987:JCP}
	{\sc Greengard, L, \& Rokhlin, V}. 1987.
	\newblock A fast algorithm for particle simulations.
	\newblock {\em J. Comput. Phys.}, {\bf 73}(2), 325--348.
	
    \bibitem[\protect\citename{FMM3DBIE_paper, }2021]{Greengard2021:JCP}
	{\sc Greengard, L., O'Neil, M., Rachh, M., \& Vico, F.} 2021.
	\newblock Fast multipole methods for the evaluation of layer potentials with locally-corrected quadratures.
	\newblock {\em J. Comput. Phys.: X}, {\bf 10}, 100092.
	
	\bibitem[\protect\citename{Gumerov \& Duraiswami, }2006]{Gumerov2006:JCP}
	{\sc Gumerov, N~A, \& Duraiswami, R}. 2006.
	\newblock Fast multipole method for the biharmonic equation in three
	dimensions.
	\newblock {\em J. Comput. Phys.}, {\bf 215}(1), 363--383.
	
	\bibitem[\protect\citename{Gumerov \& Duraiswami, }2008]{Gumerov2008:JCP}
	{\sc Gumerov, N~A, \& Duraiswami, R}. 2008.
	\newblock Fast multipole methods on graphics processors.
	\newblock {\em J. Comput. Phys.}, {\bf 227}(18), 8290--8313.
	
	\bibitem[\protect\citename{Gumerov \& Duraiswami, }2013]{Gumerov2013:JCP}
	{\sc Gumerov, N~A, \& Duraiswami, R}. 2013.
	\newblock Efficient {FMM} accelerated vortex methods in three dimensions via
	the {Lamb--Helmholtz} decomposition.
	\newblock {\em J. Comput. Phys.}, {\bf 240}(1), 310--328.
	
	\bibitem[\protect\citename{Gumerov {\em et~al.}, }2020]{Gumerov2020:IEEE}
	{\sc Gumerov, N~A, Adelman, R~N, \& Duraiswami, R}. 2020.
	\newblock Boundary Element Solution of Electromagnetic Fields for {Non-Perfect}
	Conductors at Low Frequencies and Thin Skin Depths.
	\newblock {\em IEEE Trans. Magn.}, {\bf 56}(11), 1--12.
	
	\bibitem[\protect\citename{Gumerov {\em et~al.}, }2023]{Gumerov2023:arxiv}
	{\sc Gumerov, N~A, Kaneko, S, \& Duraiswami, R}. 2023.
	\newblock Analytical Galerkin boundary integrals of Laplace kernel layer potentials in $\mathbb {R}^3$.
	\newblock {\em arXiv preprint arXiv:2302.03247 (2023)}
	
	
	\bibitem[\protect\citename{Hess \& Smith, }1967]{Hess1967:PAS}
	{\sc Hess, J~L, \& Smith, A M~O}. 1967.
	\newblock Calculation of potential flow about arbitrary bodies.
	\newblock {\em Prog. Aerosp. Sci.}, {\bf 8}, 1--138.
	
	\bibitem[\protect\citename{Kaneko, }2022]{KGD2023RIPE}
	{\sc Kaneko, S, Gumerov, N~A, \& Duraiswami, R}. 2023.
	\newblock Recursive {A}nalytical {Q}uadrature of {L}aplace and {H}elmholtz {L}ayer {P}otentials in $\mathbb{R}^3$.
	\newblock {\em arXiv preprint arXiv:2302.02196 (2023)}.
	
	\bibitem[\protect\citename{Lenoir \& Salles, }2012]{Lenoir2012:SISC}
	{\sc Lenoir, Marc, \& Salles, Nicolas}. 2012.
	\newblock Evaluation of 3-D singular and nearly singular integrals in Galerkin BEM for thin layers. \newblock{\em SIAM Journal on Scientific Computing}, {\bf 34}(6), A3057-A3078.
	
	\bibitem[\protect\citename{Newman, }1986]{Newman1986:JEM}
	{\sc Newman, John, Nicholas}. 1986.
	\newblock Distributions of sources and normal dipoles over a quadrilateral panel. \newblock{\em Journal of Engineering Mathematics}, {\bf 20}(2), 113-126.
	
	\bibitem[\protect\citename{Nishimura, }2002]{Nishimura2002:ASME}
	{\sc Nishimura, N}. 2002.
	\newblock {\em Fast multipole accelerated boundary integral equation methods}. \newblock {\em Appl. Mech. Rev.}, {\bf 55}(4), 299-324.
	
	\bibitem[\protect\citename{Nishimura {\em et~al.}, }1999]{Nishimura1999:EABE}
	{\sc Nishimura, Naoshi, Yoshida, Ken-Ichi, \& Kobayashi, Shoichi}. 1999.
	\newblock A fast multipole boundary integral equation method for crack problems
	in {3D}.
	\newblock {\em Eng. Anal. Bound. Elem.}, {\bf 23}(1), 97--105.
	
	\bibitem[\protect\citename{Of {\em et~al.}, }2005]{Of2005:CVS}
	{\sc Of, G, Steinbach, O, \& Wendland, W~L}. 2005.
	\newblock Applications of a fast multipole Galerkin in boundary element method
	in linear elastostatics.
	\newblock {\em Comput. Vis. Sci.}, {\bf 8}(3), 201--209.
	
	\bibitem[\protect\citename{Siegel \& Tornberg, }2018]{Siegel2018-xj}
	{\sc Siegel, Michael, \& Tornberg, Anna-Karin}. 2018.
	\newblock A local target specific quadrature by expansion method for evaluation
	of layer potentials in {3D}.
	\newblock {\em J. Comput. Phys.}, {\bf 364}(1), 365--392.
	
	\bibitem[\protect\citename{Tihon \& Craeye, }2018]{Tihon2018:IEEETAP}
	{\sc Tihon, Denis, \& Craeye, Christophe}. 2018.
	\newblock All-analytical evaluation of the singular integrals involved in the method of moments. \newblock{\em IEEE Transactions on Antennas and Propagation}, {\bf 66}(4), 1925-1936.

	\bibitem[\protect\citename{Vioreanu \& Rokhlin, }2014]{Vioreanu2014}
	{\sc Vioreanu, B., \& Rokhlin, V.} 2014.
	\newblock Spectra of multiplication operators as a numerical tool.
	\newblock {\em SIAM J. on Sci. Comp.}, {\bf 36}(1), A267--A288.
	
	\bibitem[\protect\citename{Wala \& Klockner, }2019]{Wala_2019-oe}
	{\sc Wala, Matt, \& Klockner, Andreas}. 2019.
	\newblock A fast algorithm for Quadrature by Expansion in three dimensions.
	\newblock {\em J. Comput. Phys.}, {\bf 388}(July), 655--689.
	
	\bibitem[\protect\citename{Wala \& Kl{\"o}ckner, }2020]{Wala2020-bb}
	{\sc Wala, Matt, \& Kl{\"o}ckner, Andreas}. 2020.
	\newblock Optimization of fast algorithms for global Quadrature by Expansion
	using target-specific expansions.
	\newblock {\em J. Comput. Phys.}, {\bf 403}(15), 108976.
	
	\bibitem[\protect\citename{Wang {\em et~al.}, }2007]{Wang2007:IJNME}
	{\sc Wang, Haitao, Lei, Ting, Li, Jin, Huang, Jingfang, \& Yao, Zhenhan}. 2007.
	\newblock A parallel fast multipole accelerated integral equation scheme for
	{3D} Stokes equations.
	\newblock {\em Int. J. Numer. Methods Eng.}, {\bf 70}(7), 812--839.
	
	\bibitem[\protect\citename{Wilton {\em et~al.}, }1984]{Wilton1984:IEEETAP}
	{\sc Wilton, D. R. S. M., Rao, S., Glisson, A. W., Schaubert, D., Al-Bundak, O., \& Butler, C.} 1984.
	\newblock Potential integrals for uniform and linear source distributions on polygonal and polyhedral domains. 
	\newblock{\em IEEE Transactions on Antennas and Propagation}, {\bf 32}(3), 276-281.
	

	
\end{thebibliography}

\end{document}